\documentclass[]{eptcs}
\usepackage{mathtools}
\usepackage{caption}
\usepackage{subcaption}
\usepackage[draft]{minted}
\usepackage[export]{adjustbox}
\usepackage{wrapfig}

\title{Operadic Modeling of Dynamical Systems: \\ Mathematics and Computation}
\author{
\hspace{0.5in}
\and
Sophie Libkind
\institute{Stanford University\\
Palo Alto, California, USA}
\email{slibkind@stanford.edu}
\and
Andrew Baas
\institute{Georgia Tech Research Institute\\
Atlanta, Georgia, USA}
\email{andrew.baas@gtri.gatech.edu}
\and
\hspace{0.5in}
\and
Evan Patterson 
\institute{Topos Institute \\ Berkeley, California, USA}
\email{evan@topos.institute}
\and
James Fairbanks
\institute{University of Florida\\ Gainesville, Florida, USA}
\email{fairbanksj@ufl.edu}
}

\def\titlerunning{Operadic Modeling of Dynamical Systems}
\def\authorrunning{ S. Libkind, A. Baas, E. Patterson, and J.P. Fairbanks}

\usepackage{fancyhdr}
\usepackage{amsmath,amssymb,amsthm}
\usepackage{tikz}
\usepackage{wrapfig} 
\usetikzlibrary{cd, arrows,automata,positioning}
\usepackage{tipa} 
\usepackage{mathtools} 
\usepackage{caption}
\usepackage{thmtools, thm-restate} 

\usetikzlibrary{decorations.pathreplacing} 

\DeclareMathOperator{\Ca}{\mathcal{C}}
\DeclareMathOperator{\Da}{\mathcal{D}}

\DeclareMathOperator{\Oa}{\mathcal{O}}

\DeclareMathOperator{\Nb}{\mathbb{N}}

\DeclareMathOperator{\Rb}{\mathbb{R}}

\DeclareMathOperator{\op}{^\text{op}}

\DeclareMathOperator{\id}{\mathsf{id}}

\DeclareMathOperator{\ob}{\mathsf{ob}}

\DeclareMathOperator{\ev}{\textbf{ev}}

\DeclareMathOperator{\Set}{\mathbf{Set}}

\DeclareMathOperator{\Finset}{\mathbf{FinSet}}
\DeclareMathOperator{\Cospan}{\mathsf{Cospan}}

\DeclareMathOperator{\Arity}{\Finset\op}

\newcommand{\Lens}[1]{\mathsf{Lens}_{#1}}
\newcommand{\Dynam}[2]{\mathsf{Dynam}^{#2}_{\mathsf{#1}}}
\newcommand{\DWD}{\CAT{DWD}}
\newcommand{\UWD}{\CAT{UWD}}
\newcommand{\CPG}{\CAT{CPG}}

\usepackage{stmaryrd}

\newcommand{\Theory}[1]{\mathsf{Th}(#1)}
\newcommand{\CAT}[1]{\mathsf{#1}}
\newcommand{\tin}{\text{in}}
\newcommand{\tout}{\text{out}}
\newcommand{\src}{\mathsf{src}}
\newcommand{\tgt}{\mathsf{tgt}}

\newcommand{\pto}{\,\cdot\kern-.1em{\to}\,}

\makeatletter
\providecommand*{\xmapstofill@}{%
  \arrowfill@{\mapstochar\relbar}\relbar\rightarrow
}
\providecommand*{\xmapsto}[2][]{%
  \ext@arrow 0395\xmapstofill@{#1}{#2}%
}

\makeatletter
\def\slashedarrowfill@#1#2#3#4#5{%
  $\m@th\thickmuskip0mu\medmuskip\thickmuskip\thinmuskip\thickmuskip
   \relax#5#1\mkern-7mu%
   \cleaders\hbox{$#5\mkern-2mu#2\mkern-2mu$}\hfill
   \mathclap{#3}\mathclap{#2}%
   \cleaders\hbox{$#5\mkern-2mu#2\mkern-2mu$}\hfill
   \mkern-7mu#4$%
}
\def\rightslashedarrowfill@{%
  \slashedarrowfill@\relbar\relbar\mapstochar\rightarrow}
\newcommand\xslashedrightarrow[2][]{%
  \ext@arrow 0055{\rightslashedarrowfill@}{#1}{#2}}
\makeatother

\tikzset{
    vert/.style={anchor=south, rotate=90, inner sep=.5mm}
} 

\newtheorem{thm}{Theorem}[section]

\theoremstyle{definition}
\newtheorem{defn}[thm]{Definition}

\newtheorem{ex}[thm]{Example}

\newtheorem{prop}[thm]{Proposition}

\newcommand{\Eul}{\mathsf{Euler}}

\newcommand{\lens}[2]{{{#1} \choose {#2}}}

\newcommand{\Euc}{\mathsf{Euc}}

\usepackage{amsthm}
\usepackage{stmaryrd}
\usepackage{newpxtext}
\usepackage[varg,bigdelims]{newpxmath}
\usepackage{tikz}

\usetikzlibrary{
	cd,
	math,
	backgrounds,
	decorations.markings,
	decorations.pathreplacing,
	positioning,
	arrows.meta,
	circuits.logic.US,
	shapes,
	calc,
	fit,
	trees,
	quotes}

\tikzset{
  WD/.style={
  	label/.style={
    	font=\everymath\expandafter{\the\everymath\scriptstyle},
      inner sep=0pt,
      node distance=2pt and -2pt},
  	label distance=-2pt,
  	every to/.style={draw},
    semithick,
    node distance=\bbx and \bby,
    decoration={markings, mark=at position \stringdecpos with \stringdec},
    bb port length=3pt,
  	bb port sep=1,
		bb inside color=white,
		bb outside color=black,
	 	bbx = .4cm,
		bb min width=.4cm,
	  bby = 2ex,
	  bb penetrate=0,
	  bb rounded corners=2pt,
	  dot size=3pt,
    shell size = 16pt,
   	penetration = 0pt,
    link size = 2pt,
    shell color = blue,
  	shell inside color=\pcolor!20,
 	  shell outside color=\pcolor!50!black,
  	surround sep=2pt,
    ar/.style={postaction={decorate}},
  	execute at begin picture={\tikzset{
  		x=\bbx, y=\bby, 
			circuit logic US, tiny circuit symbols
			}
		}
  },
  beamer/.style={
  	bbx=.4cm,
		bb min width=.4cm,
		bby=6pt,
		bb port length=3pt,
		bb port sep=.5,
  	dot size=1pt,
    shell size = 11pt, 
   	penetration = 0pt,
    link size = 1pt,
    shell color = blue,
    surround sep=1pt,
    inner sep=1pt,
    font=\tiny,
    bb inside color=\picolor,
    bb outside color=\pocolor,
	},
	bb standard colors/.style={bb inside color=white, bb outside color=black},
	bb inside color/.store in=\bbicolor,
	bb outside color/.store in=\bbocolor,
  bbx/.store in=\bbx,
  bby/.store in=\bby,
  bb port sep/.store in=\bbportsep,
  bb port length/.store in=\bbportlen,
  bb penetrate/.store in=\bbpenetrate,
  bb min width/.store in=\bbminwidth,
  bb rounded corners/.store in=\bbcorners,
  bb/.code 2 args={
    \pgfmathsetlengthmacro{\bbheight}{\bbportsep * (max(#1,#2)+1) * \bby}
    \pgfkeysalso{
      draw=\bbocolor,
      fill=\bbicolor,
      minimum height=\bbheight,
      minimum width=\bbminwidth,
      outer sep=0pt,
      rounded corners=\bbcorners,
      thick,
      prefix after command={\pgfextra{\let\fixname\tikzlastnode}},
      append after command={\pgfextra{\draw
      	\ifnum #1=0{} \else foreach \i in {1,...,#1} {
        	($(\fixname.north west)!{(2*\i-1)/(2*#1)}!(\fixname.south west)$) +(-\bbportlen,0) coordinate (\fixname_in\i) -- +(\bbpenetrate,0) coordinate (\fixname_in\i')}\fi 
        \ifnum #2=0{} \else foreach \i in {1,...,#2} {
        	($(\fixname.north east)!{(2*\i-1)/(2*#2)}!(\fixname.south east)$) +(-
\bbpenetrate,0) coordinate (\fixname_out\i') -- +(\bbportlen,0) coordinate (\fixname_out\i)}\fi
;
       }}}
		},
	dot size/.store in=\dotsize,
	dot/.style={
		circle, draw, thick, inner sep=0, fill=black, minimum width=\dotsize
	},
	bb name/.style={
    append after command={
		\pgfextra{\node[anchor=north] at (\fixname.north) {#1};}
		}
	},
  shell size/.store in=\psize,
	penetration/.store in=\penetration,
  spacing/.store in=\spacing,
  link size/.store in=\lsize,
  shell color/.store in=\pcolor,
 	shell inside color/.store in=\picolor,
 	shell outside color/.store in=\pocolor,
 	surround sep/.store in=\ssep,
 	link/.style={
  	circle, 
  	draw=black, 
  	fill=black,
  	inner sep=0pt, 
 		minimum size=\lsize
 	},
  shell/.style={
 		circle, 
 		draw = \pocolor, 
  	fill = \picolor,
  	minimum size = \psize
  },
  func/.style={
  	shell,
		rectangle,
		rounded corners=.5*\psize,
		inner ysep=.125*\psize,
		minimum width=1.125*\psize,
		inner xsep=.25*\psize,
  },
  funcr/.style={
    func,
    rectangle round north west=false, 
		rectangle round south west=false,
  },
  funcl/.style={
    func,
		rectangle round north east=false, 
		rectangle round south east=false,
  },
  funcu/.style={
    func,
		rectangle round south east=false, 
		rectangle round south west=false,
  },
  funcd/.style={
    func,
		rectangle round north east=false, 
		rectangle round north west=false,
  },
  outer shell/.style={
 		ellipse, 
 		draw,
  	inner sep=\ssep,
  	color=gray,
 	},
  intermediate shell/.style={
 		ellipse,
 		dashed, 
  	draw,
  	inner sep=\ssep,
 		color=\pocolor,
 	},
 }

\tikzset{
	oriented WD/.style={
		every to/.style={out=0,in=180,draw},
    label/.style={
    	font=\everymath\expandafter{\the\everymath\scriptstyle},
      inner sep=0pt,
      node distance=2pt and -2pt},
    semithick,
    node distance=1 and 1,
    decoration={markings, mark=at position \stringdecpos with \stringdec},
    ar/.style={postaction={decorate}},
    execute at begin picture={\tikzset{
    	x=\bbx, y=\bby,
      every fit/.style={inner xsep=\bbx, inner ysep=\bby}}}
    },
    string decoration/.store in=\stringdec,
    string decoration={\arrow{stealth};},
    string decoration pos/.store in=\stringdecpos,
    string decoration pos=.7,
    bbx/.store in=\bbx,
    bbx = 1.5cm,
    bby/.store in=\bby,
    bby = 1.ex,
    bb port sep/.store in=\bbportsep,
    bb port sep=1.5,
    bb port length/.store in=\bbportlen,
    bb port length=4pt,
    bb penetrate/.store in=\bbpenetrate,
    bb penetrate=0,
    bb min width/.store in=\bbminwidth,
    bb min width=1cm,
    bb min height/.store in=\bbminheight,
    bb min height=1cm,
    bb rounded corners/.store in=\bbcorners,
    bb rounded corners=2pt,
    bb spider/.style={
    	bb port sep=1, bb port length=10pt, bbx=.4cm, bb min width=.4cm, bby=.8ex},
    bb small/.style={
    	bb port sep=1, bb port length=2.5pt, bbx=.4cm, bb min width=.4cm, bby=.7ex},
		bb medium/.style={
			bb port sep=1, bb port length=2.5pt, bbx=.4cm, bb min width=.4cm, bby=.9ex},
    bb/.code n args={4}{
    	\pgfmathsetlengthmacro{\bbheight}{\bbportsep * (max(#1,#2)+1) * \bby}
    	\pgfmathsetlengthmacro{\bbwidth}{\bbportsep * (max(#3,#4)+1) * \bby}
      \pgfkeysalso{draw,minimum height=\bbminheight,minimum
       width=\bbminwidth,outer sep=0pt,
         rounded corners=\bbcorners,thick,
         prefix after command={\pgfextra{\let\fixname\tikzlastnode}},
         append after command={\pgfextra{\draw
            \ifnum #1=0{} \else foreach \i in {1,...,#1} {
            	($(\fixname.north west)!{\i/(#1+1)}!(\fixname.south west)$) +(-\bbportlen,0) coordinate (\fixname_in\i) -- +(\bbpenetrate, 0) coordinate (\fixname_in\i')}\fi 
            \ifnum #2=0{} \else foreach \i in {1,...,#2} {
            	($(\fixname.north east)!{\i/(#2+1)}!(\fixname.south east)$) +(-
\bbpenetrate,0) coordinate (\fixname_out\i') -- +(\bbportlen,0) coordinate (\fixname_out\i)}\fi

\ifnum #3=0{} \else foreach \i in {1,...,#3} {
            	($(\fixname.north west)!{\i/(#3+1)}!(\fixname.north east)$) +(0, \bbportlen) coordinate (\fixname_top\i) -- +(0,-\bbpenetrate) coordinate (\fixname_top\i')}\fi 
\ifnum #4=0{} \else foreach \i in {1,...,#4} {
            	($(\fixname.south west)!{\i/(#4+1)}!(\fixname.south east)$) +(0, \bbpenetrate) coordinate (\fixname_bot\i) -- +(0,-\bbportlen) coordinate (\fixname_bot\i')}\fi 
;
           }}}
		},
			bb name/.style={
     	append after command={
				\pgfextra{\node[anchor=north] at (\fixname.north) {#1};}
			}
		}
  }

  \tikzset{
  	unoriented WD/.style={
  		every to/.style={draw},
  		shorten <=-\penetration, shorten >=-\penetration,
  		label distance=-2pt,
  		thick,
  		node distance=\spacing,
  		execute at begin picture={\tikzset{
  			x=\spacing, y=\spacing}}
  		},
  	pack size/.store in=\psize,
  	pack size = 8pt,
  	spacing/.store in=\spacing,
  	spacing = 8pt,
  	link size/.store in=\lsize,
  	link size = 2pt,
		penetration/.store in=\penetration,
		penetration = 2pt,
  	pack color/.store in=\pcolor,
  	pack color = blue,
  	pack inside color/.store in=\picolor,
  	pack inside color=blue!20,
  	pack outside color/.store in=\pocolor,
  	pack outside color=blue!50!black,
  	surround sep/.store in=\ssep,
  	surround sep=8pt,
  	link/.style={
  		circle, 
  		draw=black, 
  		fill=black,
  		inner sep=0pt, 
  		minimum size=\lsize
  	},
  	pack/.style={
  		circle, 
  		draw = \pocolor, 
  		fill = \picolor,
  		inner sep = .25*\psize,
  		minimum size = \psize
  	},
  	outer pack/.style={
  		ellipse, 
  		draw,
  		inner sep=\ssep,
  		color=\pocolor,
  	},
  	intermediate pack/.style={
  		ellipse,
  		dashed, 
  		draw,
  		inner sep=\ssep,
  		color=\pocolor,
  	},
  }

\tikzset{
	spider diagram/.style={
		every to/.style={out=0, in=180, draw, thick},
		thick
	},
	dot size/.store in=\dotsize,
	dot size = 5pt,
	dot fill/.store in=\dotfill,
	dot fill = black,
	leg length/.store in=\leglen,
	leg length = 15pt,
	baby/.style={dot size = 2pt, leg length = 6pt},
	young/.style={dot size = 3pt, leg length = 10pt},
	special spider/.code n args={4}{
		\pgfkeysalso{circle, draw, thick, inner sep=0, fill=\dotfill, minimum width=\dotsize,
  		prefix after command={\pgfextra{\let\fixname\tikzlastnode}},
  		append after command={\pgfextra{
  			\ifnum #1=0{} \else {\foreach \i in {1,...,#1} {
					\tikzmath{\anglei={-90*(#1+1-2*\i)/#1};}
  				\draw [thick]
						(\fixname) .. controls 
						($(\fixname.center)-(\anglei:#3/3)$) and ($(\fixname.center)-(\anglei:#3*2/3)$) .. 
						({$(\fixname)-(\anglei:#3*2/3)$}-|{$(\fixname)-(#3,0)$}) coordinate (\fixname_in\i);
  			}}\fi
  			\ifnum #2=0{} \else {\foreach \i in {1,...,#2} {
					\tikzmath{\anglei={90*(#2+1-2*\i)/#2};}
  				\draw [thick]
						(\fixname.center) .. controls 
						($(\fixname.center)+(\anglei:#4/3)$) and ($(\fixname.center)+(\anglei:#4*2/3)$) .. 
						({$(\fixname.center)+(\anglei:#4*2/3)$}-|{$(\fixname.center)+(#4,0)$}) coordinate (\fixname_out\i);
  			}}\fi
  		}}
		}
	},
	spider/.code 2 args={
		\pgfkeysalso{special spider={#1}{#2}{\leglen}{\leglen}}
	}
}

\tikzset{Yonepart/.pic={
	\node[bb={1}{2},bb name = {\tiny$X_{11}$}] (X11) {};
	\node[bb={2}{2},below right=of X11,bb name = {\tiny$X_{12}$}] (X12) {};
	\node[bb={2}{1}, above right=of X12,bb name = {\tiny$X_{13}$}] (X13) {};
	\node[bb={2}{2}, fit={($(X11.north west)+(.3,1.5)$) (X12)  ($(X13.east)+(-.3,0)$)},bb name = {\scriptsize $Y_1$}] (Y1) {};
	\draw (Y1_in1') to (X11_in1);	
	\draw (Y1_in2') to (X12_in2);
	\draw (X11_out1) to (X13_in1);
	\draw (X11_out2) to (X12_in1);
	\draw (X12_out1) to (X13_in2);
	\draw (X12_out2) to (Y1_out2');
	\draw (X13_out1) to (Y1_out1');
	\coordinate (bottombox) at ($(X12.south)$);
	\coordinate (rightbox) at ($(X13.east)$);
	\coordinate (Y1northwest) at ($(Y1.north west)$);
	}
}

\tikzset{Ytwopart/.pic={
	\node[bb={2}{2}, bb name = {\tiny$X_{21}$}] (X21) {};
	\node[bb={1}{2},above right=-1 and 1 of X21,bb name = {\tiny$X_{22}$}] (X22) {};
	\node[bb={1}{2}, fit={($(X21.south west)+(-.25,0)$) ($(X22.north east)+(.25,3.5)$)},bb name = {\scriptsize$Y_2$}] (Y2){};
	\draw (Y2_in1') to (X21_in2);
	\draw (X21_out1) to (X22_in1);
	\draw (X22_out2) to (Y2_out1');
	\draw let \p1=(X22.south east), \p2=($(Y2_out2)$), \n1={\y1-\bby}, \n2=\bbportlen in
	  (X21_out2) to (\x1+\n2,\n1) -- (\x1+\n2,\n1) to (Y2_out2');
	\draw let \p1=(X22.north east), \p2=(X21.north west), \n1={\y1+\bby}, \n2=\bbportlen in
          (X22_out1) to[in=0] (\x1+\n2,\n1) -- (\x2-\n2,\n1) to[out=180] (X21_in1);
          }
}

\tikzset{SmallNeuronPic/.pic={
 \node[bb={3}{1}] (N1) {$\scriptstyle N_1$};
  \node[bb={2}{1}, above right=.7 and 3.5 of N1] (N2) {$\scriptstyle N_2$};
  \node[bb={2}{1}, below =of N2] (N3) {$\scriptstyle N_3$};
  \node[bb={3}{1}, below =of N3] (N4) {$\scriptstyle N_4$};
  \node[bb={6}{8}, fit={($(N1.west)-(.5,0)$) ($(N2.north)+(0,2)$) ($(N3.east)+(1.5,0)$) ($(N4.south)-(0,1)$)}, bb name={$\scriptstyle X$}] (X) {};
  \draw (X_in1') to (N2_in1);
  \draw (X_in2') to (N1_in1);
  \draw (X_in3') to (N1_in2);
  \draw (X_in4') to (N1_in3);
  \draw (X_in6') to (N4_in2);
  \draw (N1_out1) to (N2_in2);
  \draw (N1_out1) to (N3_in1);
  \draw (N1_out1) to (N4_in1);
  \draw (N2_out1) to (X_out1');
  \draw (N2_out1) to (X_out2');
  \draw (N2_out1) to (X_out3');
  \draw (N3_out1) to (X_out4');
  \draw (N3_out1) to (X_out5');
  \draw (N3_out1) to (X_out6');
  \draw (N4_out1) to (X_out7');
  \draw (N4_out1) to (X_out8'); 
  \draw (X_in5') to[looseness=2] (N3_in2);
  \draw let \p1=(N4.south east), \p2=(N4.south west), \n1={\y2-\bby}, \n2=\bbportlen in
          (N3_out1) to[in=0] (\x1+\n2,\n1) -- (\x2-\n2,\n1) to[out=180] (N4_in3);
}
}

\tikzset{SmallNeuronDashed/.pic={
 \node[bb={3}{1}] (N1) {$\scriptstyle N_1$};
  \node[bb={2}{1}, above right=.7 and 3.5 of N1] (N2) {$\scriptstyle N_2$};
  \node[bb={2}{1}, below =of N2] (N3) {$\scriptstyle N_3$};
  \node[bb={3}{1}, below =of N3] (N4) {$\scriptstyle N_4$};
  \node[bb={6}{8}, fit={($(N1.west)-(.5,0)$) ($(N2.north)+(0,2)$) ($(N3.east)+(1.5,0)$) ($(N4.south)-(0,1)$)}, bb name={$\scriptstyle X$}] (X) {};
  \draw[dashed] (X_in1') to (N2_in1);
  \draw[dashed] (X_in2') to (N1_in1);
  \draw[dashed] (X_in3') to (N1_in2);
  \draw[dashed] (X_in4') to (N1_in3);
  \draw[dashed] (X_in6') to (N4_in2);
  \draw[dashed] (N1_out1) to (N2_in2);
  \draw[dashed] (N1_out1) to (N3_in1);
  \draw[dashed] (N1_out1) to (N4_in1);
  \draw[dashed] (N2_out1) to (X_out1');
  \draw[dashed] (N2_out1) to (X_out2');
  \draw[dashed] (N2_out1) to (X_out3');
  \draw[dashed] (N3_out1) to (X_out4');
  \draw[dashed] (N3_out1) to (X_out5');
  \draw[dashed] (N3_out1) to (X_out6');
  \draw[dashed] (N4_out1) to (X_out7');
  \draw[dashed] (N4_out1) to (X_out8'); 
  \draw[dashed] (X_in5') to[looseness=2] (N3_in2);
  \draw[dashed] let \p1=(N4.south east), \p2=(N4.south west), \n1={\y2-\bby}, \n2=\bbportlen in
          (N3_out1) to[in=0] (\x1+\n2,\n1) -- (\x2-\n2,\n1) to[out=180] (N4_in3);
}
}

\tikzset{SmallNestingPic/.pic={
\path (0,0) pic [purple] {Yonepart};
\path ($(rightbox)+(5,-5)$) pic [orange] {Ytwopart};
 
\node[bb={1}{2}, fit={($(Y1northwest)+(-.5,4)$) ($(Y2.south east)+(1,0)$)}, bb name={\small $Z$}] (Z) {};
\draw (Z_in1') to (Y1_in2);
\draw let \p1=(Y2.north west),\p2=(Y2.north east),\n1={\y2+\bby},\n2=\bbportlen in
  (Y1_out1) to (\x1+\n2,\n1)--(\x2+\n2,\n1) to (Z_out1');
\draw (Y1_out2) to (Y2_in1);
\draw (Y2_out2) to (Z_out2');
\draw let \p1=(Y2.north east), \p2=(Y1.north west), \n1={\y2+\bby}, \n2=\bbportlen in
          (Y2_out1) to[in=0] (\x1+\n2,\n1) -- (\x2-\n2,\n1) to[out=180] (Y1_in1);
          }
}

\tikzset{Zredgreen/.pic={
\node[bb={2}{2}, green!50!black, bb name = $\scriptstyle Y_1$] (YY1) {};
\node[bb={1}{2}, red, below right=-1 and 2 of YY1, bb name=$\scriptstyle Y_2$] (YY2) {};
\node[bb={1}{2}, fit={($(YY1.north west)+(-.5,4)$) ($(YY2.south east)+(.5,-2)$)}, bb name={\scriptsize $Z$}] (Z) {};
\draw (Z_in1') to (YY1_in2);
\draw (YY1_out1) to (Z_out1');
\draw (YY1_out2) to (YY2_in1);
\draw (YY2_out2) to (Z_out2');
\draw let \p1=(YY2.north east), \p2=(YY1.north west), \n1={\y2+\bby}, \n2=\bbportlen in
          (YY2_out1) to[in=0] (\x1+\n2,\n1) -- (\x2-\n2,\n1) to[out=180] (YY1_in1);
}
}

\tikzset{Zcombined/.pic={
	\node[bb={1}{2},green!25!black,bb name = {\tiny$X_{11}$}] (X11) {};
	\node[bb={2}{2},green!25!black,below right=of X11,bb name = {\tiny$X_{12}$}] (X12) {};
	\node[bb={2}{1}, green!25!black,above right=of X12,bb name = {\tiny$X_{13}$}] (X13) {};
	\draw (X11_out1) to (X13_in1);
	\draw (X11_out2) to (X12_in1);
	\draw (X12_out1) to (X13_in2);

	\node[bb={2}{2}, red!30!black, below right = 0 and 1.25 of X12, bb name = {\tiny$X_{21}$}] (X21) {};
	\node[bb={1}{2}, red!30!black, above right=-1 and 1 of X21,bb name = {\tiny$X_{22}$}] (X22) {};
	\draw (X21_out1) to (X22_in1);
	\draw let \p1=(X22.north east), \p2=(X21.north west), \n1={\y1+\bby}, \n2=\bbportlen in
          (X22_out1) to[in=0] (\x1+\n2,\n1) -- (\x2-\n2,\n1) to[out=180] (X21_in1);
        
        \node[bb={1}{2}, fit = {($(X11.north east)+(-1,3)$) (X12) (X13) ($(X21.south)+(0,-1)$) ($(X22.east)+(.5,0)$)}, bb name ={\scriptsize $Z$}] (Z) {};
	
	\draw (Z_in1') to (X12_in2);
	\draw (X13_out1) to (Z_out1');
	\draw (X12_out2) to (X21_in2);
	\draw let \p1=(X22.south east),\n1={\y1-\bby}, \n2=\bbportlen in
	  (X21_out2) to (\x1+\n2,\n1) to (Z_out2');
	\draw let \p1=(X22.north east), \p2=(X11.north west), \n1={\y2+\bby}, \n2=\bbportlen in
          (X22_out2) to[in=0] (\x1+\n2,\n1) -- (\x2-\n2,\n1) to[out=180] (X11_in1);
}
}

\begin{document}
\maketitle

\begin{abstract}
Dynamical systems are ubiquitous in science and engineering as models of phenomena that evolve over time. Although complex dynamical systems tend to have important modular structure, conventional modeling approaches suppress this structure. Building on recent work in applied category theory, we show how deterministic dynamical systems, discrete and continuous, can be composed in a hierarchical style. In mathematical terms, we reformulate some existing operads of wiring diagrams and introduce new ones, using the general formalism of $\Ca$-sets (copresheaves). We then establish dynamical systems as algebras of these operads. In a computational vein, we show that Euler's method is functorial for undirected systems, extending a previous result for directed systems. All of the ideas in this paper are implemented as practical software using Catlab and the AlgebraicJulia ecosystem, written in the Julia programming language for scientific computing.
\end{abstract}

\section{Introduction}

Category theory is about finding the right abstractions---identifying the salient, general features of the objects of study. In applied category theory the chosen objects of study lie outside of pure mathematics. 
One important thread of finding the right abstractions in the sciences has been understanding the composition of dynamical systems. Dynamical systems are a general and ubiquitous class of models which capture changing phenomena. For example, automata model the changing of states in a computer, Petri nets model the changing concentrations of chemicals in a reaction network, and flows on a manifold model the evolution of physical systems. 
There is a long history of scientists developing ad hoc graphical languages for specifying, communicating, and refining large composite models~\cite{hall1977ecosystem}.
However, these languages are informal, and so the modular structure of a complex system is often lost in implementation and cannot be used for model calibration or analysis. 
In this paper, we demonstrate how operads and operad algebras can be used to formalize compositional modeling, and we present a Julia package for dynamical systems that preserves the compositional structure.
If applied category theory is about \textit{finding} the right abstractions for science, then the present work exemplifies \textit{implementing} the right abstractions for science.

Existing work on composing dynamical systems varies along two axes. The first axis is semantic: \textit{what} is a dynamical system? Dynamical systems are an extremely broad class of models, and previous work falls on many different points along the semantic axis. These points include circuit diagrams, Petri nets, Markov processes, finite state automata, ODEs, hybrid systems, and Lagrangian and Hamiltonian systems~\cite{fong2015DecoratedCospans, baez2017CompositionalFramework,baez2016compositionalframework,vagner2015AlgebrasOpen, libkind2020AlgebraResource, lerman2020Networkshybrid, baez2021Opensystems}. In this paper we focus on two kinds of dynamics:  continuous flows and discrete transitions.

The second axis is syntactic: \textit{how} do dynamical systems compose? Two distinct styles of composition have emerged: directed and undirected, also called machine composition and resource sharing. In directed composition, information is transferred from designated senders to designated receivers. Systems are driven by the behavior of other systems but otherwise have independent dynamics.  In undirected composition, systems compose by sharing resources or observations. Composed systems affect each other only by acting on the shared medium. An important distinction is that undirected composition is not equivalent to symmetric directed composition. The directed and undirected perspectives are unified in \cite{libkind2020AlgebraResource}.

\subsection{Contributions}

\begin{enumerate}
    \item A practical implementation of operads and their algebras in the programming language Julia,\footnote{The software implementation can be found at \url{https://github.com/AlgebraicJulia/AlgebraicDynamics.jl}.} which is widely used for scientific computing.
    
    \item A reformulation of previously studied operads using $\Ca$-sets, a diagrammatic approach to defining data structures. A new instance of this abstraction, the operad of circular port graphs, is also introduced.
    
    \item Two new algebras for composing dynamical systems. The first algebra represents a directed composition of continuous and discrete systems that extends the syntax of algebras previously studied in~\cite{spivak2016steadystates, schultz2020DynamicalSystems, vagner2015AlgebrasOpen} to include merging and creating wires. The second algebra represents the undirected composition of discrete dynamical systems. 
    
    \item A proof that Euler's method is functorial for undirected systems, plus an implementation in Julia of functorial Euler's method for both directed and undirected systems.

\end{enumerate}

\noindent {\bf Acknowledgements} The authors were supported by DARPA Awards W911NF2010292 and HR00112090067 along with AFOSR Award FA9550-20-1-0348. The authors thank Micah Halter and Owen Lynch for support in developing Catlab and the AlgebraicDynamics software packages. They also thank David Spivak for his helpful insights into operads of wiring diagrams.

\section{Preliminary Definitions}
\subsection{Operads and Operad Algebras}
Operads and operad algebras formalize notions of syntax and semantics. In contrast, modeling tools  generally obscure the distinction between syntax and semantics. These blurred lines make it challenging to interoperate between modeling frameworks and to independently adjust model syntax and model semantics. In this section, we give the mathematical background for operads and operad algebras which  form the foundation of our software implementation.

Throughout we use \textit{operads} to refer to symmetric colored operads or equivalently symmetric multicategories. We will also refer to the objects of an operad as its \textit{types} and the morphisms of an operad as its \textit{terms} in order to highlight the connection with syntax.

\begin{defn}
An \emph{operad} $\Oa$ consists of a collection of types $\ob \Oa$ and for each $n \in \Nb^+$ and types $s_1,\dots, s_n, t \in \ob \Oa$, a collection of terms $\Oa(s_1,\dots, s_n ; t)$, along with
\begin{itemize}
    \item an identity term $1_t \in \Oa(t;t)$ for each type $t \in \ob \Oa$,
    \item substitution maps $$\circ_i: \Oa(r_1,\dots, r_m; s_i) \times \Oa(s_1,\dots, s_n;t) \to \Oa(s_1,\dots, s_{i -1}, r_1,\dots, r_m, s_{i + 1},\dots, s_n; t),$$
\end{itemize}
and permutation maps satisfying associativity, unitality, and symmetry laws.

An \emph{operad functor} $F: \Oa \to \Oa'$ is a map on types and on terms that commutes with the identity, substitution,  and permutation. Operads and their functors form a category  $\mathsf{Oprd}$.\footnote{See~\cite{spivak2013operadwiring} for a detailed exposition of operads and operad functors that aligns with their usage here. See also~\cite{leinster2004Higheroperads}.}
\end{defn}

A rich source of operads and operad algebras is the category $\mathsf{SMC}$ of symmetric monoidal categories (SMCs) and lax monoidal functors. Specifically, there exists a functor $\Oa: \mathsf{SMC} \to \mathsf{Oprd}$ sending each SMC $(\Ca, \otimes, 1)$ to its underlying operad $\Oa(\Ca)$ with types $\ob \Ca$ and terms $\Oa(\Ca)(s_1, \dots, s_m; t) := \Ca(s_1 \otimes \cdots \otimes s_m, t)$~\cite{leinster2004Higheroperads}. When clear from context, we denote $\Oa(\Ca)$ simply by $\Ca$.  

\begin{defn}
Given an operad $\Oa$, an \emph{algebra of $\Oa$} or simply an \emph{$\Oa$-algebra} is an operad functor $F: \Oa \to \Oa(\Set)$. We call a pair $(t \in \ob \Oa, m \in Ft) \in \int_{\ob \Oa} F$ an \emph{element} of the algebra.
\end{defn}

Symmetric monoidal categories ease the way for mathematical formalization and analysis. In this work, all of the operads and operad algebras for modeling dynamical systems are induced by symmetric monoidal categories and lax monoidal functors, respectively. However, the operadic perspective is better suited to computing because it directly supports $n$-ary operations rather than requiring that $n$-ary operations be decomposed into a tree of binary operations.  The operadic viewpoint is highlighted in the Julia implementation (Section~\ref{sec:implementation}). Now we give two SMCs whose underlying operads define syntaxes for directed and undirected composition of dynamical systems, respectively.

\begin{figure}[htbp]
    \centering
    \begin{subfigure}[b]{0.15\textwidth}
    \centering
    \begin{tikzpicture}[oriented WD, bbx = 1cm, bby =.5cm, bb port length=4pt, bb port sep=.75, bb min width=1cm, bb min height=1cm]
    \node[bb={0}{0}{2}{3}] (B){};

\end{tikzpicture}
    \subcaption{}
    \end{subfigure}\hfill
    \begin{subfigure}[b]{0.4\textwidth}
    \centering
    \begin{tikzpicture}[oriented WD, bbx = .5cm, bby =.25cm, bb port length=4pt, bb port sep=.75, bb min width=1cm, bb min height=1cm]
    \node[bb={0}{0}{2}{0}] (A){};
    \node[bb={0}{0}{1}{2}, right=of A] (B){};
    \node[bb={0}{0}{0}{3}, right=of B] (C){};
    
    \node[bb={0}{0}{2}{2}, fit={($(A.north west)+(0,2)$)($(C.south east)+(0, -3.5)$)}] (tot){};
    
    \draw[ar, color=violet] (tot_top1') to [out=-90, in=90] (A_top2);

    \draw[ar, color=violet] let \p1=(A.south west), \p2=(A.north west), \n1=\bbportlen, \n2=\bby in
     (B_bot1') to [out=-90, in=-90] (\x1 - \bby, \y1 - \n1) to [out=90, in=-90] (\x2 - \bby, \y2 +\n1) to [out=90, in=90] (A_top1);
    
    \draw[ar, color=violet] let \p1=(B.south east), \p2=(B.north east), \n1=\bbportlen, \n2=\bby in
     (B_bot2') to [out=-90, in=-90] (\x1 + \bby, \y1 - \n1) to [out=90, in=-90] (\x2 + \bby, \y2 +\n1) to [out=90, in=90] (B_top1);
     
    \draw[ar, color=orange] (B_bot2') to [out=-90, in=90] (tot_bot2);
    \draw[ar, color=orange] (C_bot2') to [out=-90, in=90] (tot_bot1);

    \node[above=2pt of tot_top1'] () {};
    \node[above=2pt of tot_top2'] () {};
    
    \node[below=2pt of A_top1'] () {};
    \node[below=2pt of A_top2'] () {};
     \node[below=2pt of B_top1'] () {};

\end{tikzpicture}
    \subcaption{}
    \end{subfigure} \hfill
    \begin{subfigure}[b]{0.4\textwidth}
    \centering
    \begin{tikzpicture}[oriented WD, bbx = .5cm, bby =.25cm, bb port length=4pt, bb port sep=.75, bb min width=1cm, bb min height=1cm]
    \node[bb={0}{0}{2}{0}] (A){};
    \node[bb={0}{0}{2}{2}, right=of A] (B){};
    \node[bb={0}{0}{0}{3}, right=of B] (C){};
    
    \node[bb={0}{0}{2}{2}, fit={($(A.north west)+(0,2)$)($(C.south east)+(0, -3.5)$)}] (tot){};
    
    \draw[ar, color=violet] (tot_top1') to [out=-90, in=90] (A_top2);
    \draw[ar, color=violet] (tot_top1') to [out=-90, in=90] (B_top1);
    
    \draw[ar, color=violet] let \p1=(A.south west), \p2=(A.north west), \n1=\bbportlen, \n2=\bby in
     (B_bot1') to [out=-90, in=-90] (\x1 - \bby, \y1 - \n1) to [out=90, in=-90] (\x2 - \bby, \y2 +\n1) to [out=90, in=90] (A_top1);
    
    \draw[ar, color=violet] let \p1=(B.south east), \p2=(B.north east), \n1=\bbportlen, \n2=\bby in
     (B_bot2') to [out=-90, in=-90] (\x1 + \bby, \y1 - \n1) to [out=90, in=-90] (\x2 + \bby, \y2 +\n1) to [out=90, in=90] (B_top1);
     
    \draw[ar, color=orange] (B_bot2') to [out=-90, in=90] (tot_bot2);
    \draw[ar, color=orange] (C_bot2') to [out=-90, in=90] (tot_bot1);
    \draw[ar, color=orange] (C_bot3') to [out=-90, in=90] (tot_bot2);
    
    \node[above=2pt of tot_top1'] () {};
    \node[above=2pt of tot_top2'] () {};
    
    \node[below=2pt of A_top1'] () {};
    \node[below=2pt of A_top2'] () {};
     \node[below=2pt of B_top1'] () {};

\end{tikzpicture}
    \subcaption{}
    \end{subfigure} 

    \caption{The graphical representation of $\DWD$ and its wide suboperad $\Oa(\Lens{\Arity})$.  (a) The type $\lens{2}{3}$. (b) A term $\lens{2}{0} + \lens{1}{2} + \lens{0}{3} \leftrightarrows \lens{2}{2}$ in $\Oa(\Lens{\Arity})$. (c) A morphism $\lens{f^\#}{f}: \lens{2}{0} + \lens{2}{2} + \lens{0}{3} \leftrightarrows \lens{2}{2}$ in $\DWD$. The orange and purple wires represent the apexes of $f$ and $f^\#$ respectively. This syntactic diagram depicts the merging of wires (e.g., the first in-port of the second inner box) and the creation of wires (e.g., the second in-port of the second inner box).  In contrast to the interpretation of a string diagram, here boxes represent types (objects) and the assembly of wires represents a term (morphism).}
    \label{fig:operad-dwd}
\end{figure}
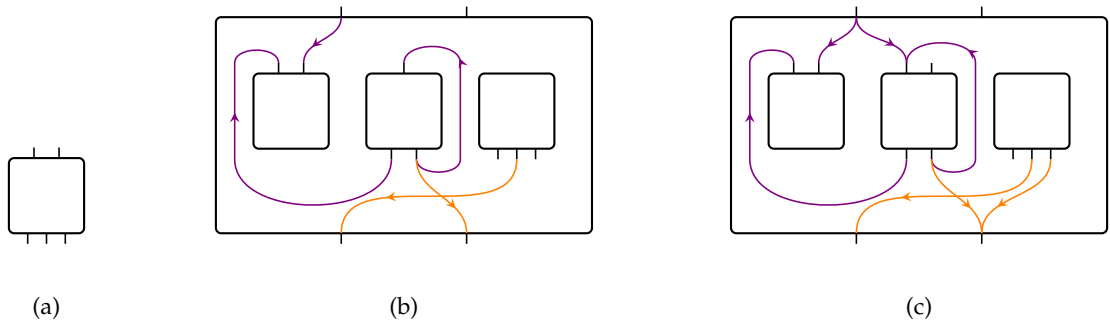

\begin{ex}[Operad of directed wiring diagrams]\label{ex:operad-dwd}
Given a category $\Ca$ with finite products, there is a lens category, $\Lens{\Ca}$, whose objects are pairs $\lens{X_\tin}{X_\tout}$ where $X_\tin, X_\tout \in \ob \Ca$ and whose morphisms $\lens{f}{f^\#} \in \Lens{\Ca}\left( \lens{X_\tin}{X_\tout},  \lens{Y_\tin}{Y_\tout} \right)$ are pair of morphisms $f: X_\tout \to Y_\tout$, $f^\#: Y_\tin \times X_\tout \to X_\tin$ ~\cite[Definition 2.2]{spivak2020GeneralizedLens}. The cartesian monoidal structure on $\Ca$ induces a symmetric monoidal structure on $\Lens{\Ca}$ \cite{moeller2020MonoidalGrothendieck}. Therefore $\Lens{\Ca}$ has an underlying operad.

In~\cite{vagner2015AlgebrasOpen, schultz2020DynamicalSystems, jazmyers2021DoubleCategories} the directed syntax for composing dynamical systems is defined by the operad underlying $\left(\Lens{\Arity}, +, \lens{0}{0}\right)$, often referred to as the operad of wiring diagrams. Following the Catlab implementation, we instead focus on the operad underlying $\left(\Lens{\Cospan(\Arity)}, + , \lens{0}{0}\right)$. We define $\DWD \coloneqq \Oa(\Lens{\Cospan(\Arity)})$ and call it the \textit{operad of directed wiring diagrams}. 
In contrast to $\Oa(\Lens{\Arity})$, which allows only copying and deletion of wires, the syntax defined by $\DWD$ can also represent merging and creation of wires.

The graphical representation of $\DWD$ extends the standard graphical representation of $\Oa(\Lens{\Arity})$ (Figure~\ref{fig:operad-dwd}). Types $\lens{X_\tin}{X_\tout}$ are represented as boxes with in-ports $X_\tin$ and out-ports $X_\tout$. Let $\lens{f^\#}{f}: \lens{X_\tin}{X_\tout} \leftrightarrows \lens{Y_\tin}{Y_\tout}$ be a term in $\DWD$. The morphism $f = X_\tout \leftarrow V \to Y_\tout$ represents a set $V$ of wires with sources and targets given by the left and right legs of the span. Likewise, $f^\#: Y_\tin + X_\tout \leftarrow W \to X_\tin$ represents another set $W$ of wires. The graphical representation emphasizes the operadic structure by having a separate box for each type in the term's domain.

\end{ex}

\begin{ex}[Operad of undirected wiring diagrams]\label{ex:operad-uwd}

We define $\UWD$ to be the operad underlying the symmetric monoidal category $(\Cospan(\Finset), + , 0)$ and call it the \textit{operad of undirected wiring diagrams}. Graphically, a type $M$ is represented by a box with $M$ exposed ports and a term $M \to J \leftarrow N$ is represented by $J$ junction nodes with wires connecting ports $M$ and $N$ to junctions according to the legs of the cospan ~\cite{spivak2013operadwiring}.

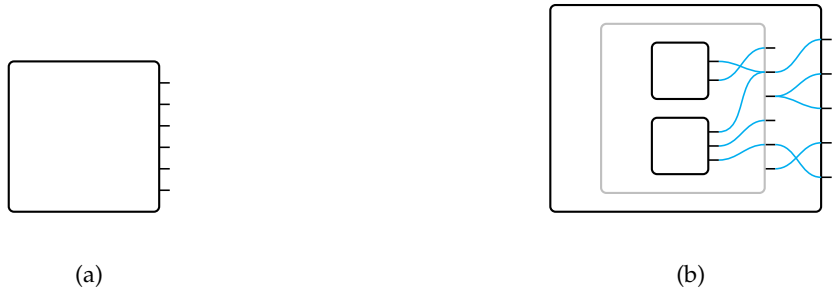
\begin{figure}[htbp]
    \centering
    \begin{subfigure}[b]{0.5\textwidth}
    \centering
    \begin{tikzpicture}[ oriented WD, bbx = 1cm, bby =.5cm, bb port length=4pt, bb port sep=.75, bb min width=2cm, bb min height = 2cm]
    \node[bb={0}{6}{0}{0}] (){};
\end{tikzpicture}
    \subcaption{}
    \end{subfigure}\hfill
    \begin{subfigure}[b]{0.5\textwidth}
    \centering
    \begin{tikzpicture}[ oriented WD, bbx =.75cm, bby =.25cm, bb port length=4pt, bb port sep=.75, bb min width = .75cm, bb min height = .75cm]
    \node[bb={0}{2}{0}{0}] (b1){};
    \node[bb={0}{3}{0}{0}, below= of b1] (b2){};
    \node[bb={0}{6}{0}{0}, fit={($(b1.north west)+(.1, 0)$)(b2.south east)}, color=lightgray] (q){};
    \node[bb={0}{5}{0}{0}, fit={($(q.north west)+(.1, 0)$)(q.south east)}] (m){};

    \draw[color=cyan] (b1_out1) edge [out=0, in=180] (q_out2');
    \draw[color=cyan] (b1_out2) edge [out=0, in=180] (q_out1');
    \draw[color=cyan] (b2_out1) edge [out=0, in=180] (q_out2');
    \draw[color=cyan] (b2_out2) edge [out=0, in=180] (q_out4');
    \draw[color=cyan] (b2_out3) edge [out=0, in=180] (q_out5');
    
    \draw[color=cyan] (m_out1') edge [out=180, in=0] (q_out2);
    \draw[color=cyan] (m_out2') edge [out=180, in=0] (q_out3);
    \draw[color=cyan] (m_out3') edge [out=180, in=0] (q_out3);
    \draw[color=cyan] (m_out4') edge [out=180, in=0] (q_out6);
    \draw[color=cyan] (m_out5') edge [out=180, in=0] (q_out5);

\end{tikzpicture}
    \subcaption{}
    \end{subfigure}

    \caption{The graphical representation of $\UWD$. (a) The type $6$. (b) A term $2 + 3 \to 6 \leftarrow 5 \in \UWD(2,3; 5)$.}
    \label{fig:operad-uwd}
\end{figure}
\end{ex}

\subsection{\texorpdfstring{$\Ca$}{C}-sets}
$\Ca$-sets are a powerful abstraction for capturing data of a fixed shape \cite{reyes2004GenericFigures,schultz2017AlgebraicDatabases}. In this section, we define $\Ca$-sets and introduce specific $\Ca$-sets implementing operad terms. 

\begin{defn}
Let $\Ca$ be a small category. An \emph{$\Ca$-set} or \emph{instance of $\Ca$} is a copresheaf over $\Ca$, equivalently a functor $X: \Ca \to \Set$. If $X$ factors through $\Finset$, then we say that $X$ is \emph{finite}.
\end{defn}

\begin{defn}
Let $F: \Ca \to \Da$ be a functor. Then there is a \emph{pullback data migration functor} $\Delta_F: [\Da, \Set] \to [\Ca, \Set]$ given by precomposition with $F$. The action of $\Delta_F$ on objects turns instances of $\Da$ into instances of $\Ca$.
\end{defn}

The category $\Ca$ is a schema that structures data, and an instance $X: \Ca \to \Set$ is an instance of the data structure.   Next, we present schemata for undirected wiring diagrams, directed wiring diagrams, and circular port graphs. In our examples, the schemata are  finitely presented categories and the finite instances of each schema comprise the terms of an operad.  Thus, we can take advantage of the rich mathematical structure of $\Ca$-sets, such as functorial data migration and the existence of finite limits and colimits, to build syntactic terms. 

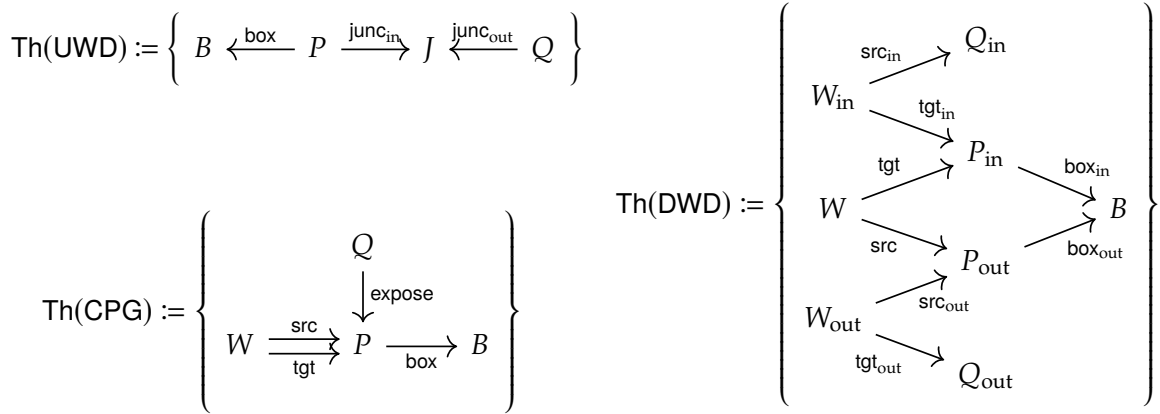
\begin{figure}
    \begin{minipage}{0.45\linewidth}
        \begin{equation*}
\Theory{\CAT{UWD}} \coloneqq \left\{
    \begin{tikzcd}
    B &P \arrow[l,swap, "\mathsf{box}"] \arrow[r, "\mathsf{junc}_\tin"] & J &  \arrow[l, "\mathsf{junc}_\tout", swap] Q
    \end{tikzcd}
  \right\}
\end{equation*}
        \begin{equation*}
\Theory{\CAT{CPG}} \coloneqq \left\{
    \begin{tikzcd}
    & Q \arrow[d, "\mathsf{expose}"] &\\
      W \arrow[r, shift left=1, "\src"]
        \arrow[r, swap, shift right=1, "\tgt"]
      & P \arrow[r, "\mathsf{box}", swap] &B \\
    \end{tikzcd}
  \right\}
 \end{equation*}
    \end{minipage} \begin{minipage}{0.55\linewidth}
        \begin{equation*}
\Theory{\CAT{DWD}} \coloneqq \left\{
    \begin{tikzcd}[row sep = .05in]
     & Q_\tin&\\
     W_\tin  \arrow[ur, "\src_\tin"] \arrow[dr, "\tgt_\tin"]&& \\
     & P_\tin \arrow[dr, "\mathsf{box}_\tin"]&\\
     W \arrow[ur,"\tgt"] \arrow[dr, "\src", swap] &&B \\
     & P_\tout \arrow[ur, "\mathsf{box}_\tout", swap]&\\
     W_\tout \arrow[ur, "\src_\tout",swap] \arrow[dr, "\tgt_\tout", swap] &&\\
     & Q_\tout &\\
    \end{tikzcd}
  \right\}
 \end{equation*}
    \end{minipage}
    \caption{The schemata for the theories of undirected wiring diagram, directed wiring diagrams, and  circular port graphs. }
    \label{fig:theories}
\end{figure}

\begin{ex}[Theory of undirected wiring diagrams]\label{ex:theory-uwd}
The schema for undirected wiring diagrams is $\Theory{\CAT{UWD}}$, defined  in Figure~\ref{fig:theories}. An instance $X$ of $\Theory{\CAT{UWD}}$ consists of a set of boxes $XB$, ports $XP$, outer ports $XQ$, and junctions $XJ$. Each box $b \in XB$ has ports $X\mathsf{box}^{-1}(b) \subseteq XP$. Each port $p \in XP$ connects to the junction $X\mathsf{junc}_\tin(p)$, and likewise for outer ports. A finite instance $X$ of $\Theory{\UWD}$ presents a term in $\Oa(\Cospan(\Finset))$ that has domain types $XB$ and underlies the morphism $XP \xrightarrow{X\mathsf{junc}_\tout} XJ \xleftarrow{X\mathsf{junc}_\tout} XQ$ in $\Cospan(\Finset)$.  Up to relabeling of the box elements and permutation of the domain types, finite instances of $\Theory{\UWD}$ correspond one-to-one with  terms of $\UWD$.
\end{ex}

\begin{ex}[Theory of directed wiring diagrams]\label{ex:theory-dwd}

The schema for directed wiring diagrams is $\Theory{\CAT{DWD}}$, defined  in Figure~\ref{fig:theories}. An instance $X$ of $\Theory{\CAT{DWD}}$ consists of a set of boxes $XB$, sets of inner in-ports and out-ports $XP_\tin$ and $XP_\tout$, sets of outer in-ports and out-ports $XQ_\tin$ and $XQ_\tout$, and a set of wires $XW_\tin + XW + XW_\tout$. Each wire has source and target given by $X\src_\tin + X\src + X\src_\tout$ and $X\tgt_\tin + X\tgt + X\tgt_\tout$ respectively. A finite instance $X$ of $\Theory{\DWD}$ presents a term in $\DWD$ that has domain types $XB$ and underlies the morphism 
$$\lens{XQ_\tin + XP_\tout \xleftarrow{X\src_\tin + X\src} XW_\tin + XW \xrightarrow{(X\tgt_\tin, X\tgt)} XP_\tin}{XP_\tout \xleftarrow{X\src_\tout} XW_\tout \xrightarrow{X\tgt_\tout} XQ_\tout} :\lens{XP_\tin}{XP_\tout} \leftrightarrows \lens{XQ_\tin}{XQ_\tout}$$ in $\Lens{\Cospan(\Arity)}$. Up to relabeling of the box elements and permutation of the domain types, finite instances of $\Theory{\DWD}$ correspond one-to-one with terms of $\DWD$.
\end{ex}

\begin{ex}[Theory of circular port graphs]\label{ex:theory-cpg} 
The schema for circular port graphs is $\Theory{\CAT{CPG}}$, defined in Figure~\ref{fig:theories}. An instance $X$  of $\Theory{\CAT{CPG}}$ consists of a set of ports $XP$, outer ports $XQ$, and wires $XW$ whose source and target are specified by $X\src$ and $X\tgt$ respectively. 

Every circular port graph induces a directed wiring diagram by functorial data migration. Let the functor $F: \Theory{\CAT{DWD}} \to \Theory{\CAT{CPG}}$ be defined on objects by $Q_\tin, Q_\tout, W_\tin, W_\tout \mapsto Q$,  $P_\tin, P_\tout \mapsto P$, $W \mapsto W$, and $B \mapsto B$  and on morphisms by $\src_\tin, \tgt_\tout \mapsto \id_{Q}$,  $\tgt_\tin, \src_\tout \mapsto \mathsf{expose}$,  $\src \mapsto \src$,  $\tgt \mapsto \tgt$, and $\mathsf{box}_\tin, \mathsf{box}_\tout \mapsto \mathsf{box}$. The pullback data migration functor $\Delta_F:  [\Theory{\CAT{CPG}}, \Set] \to [\Theory{\CAT{DWD}}, \Set]$ interprets circular port graphs as directed wiring diagrams by duplicating every port with one copy interpreted as an in-port and the other as an out-port. Example~\ref{ex:theory-dwd} gives a correspondence between finite instances of $\Theory{\DWD}$ and terms of $\DWD$. Composing $\Delta_F$ with this correspondence defines a map from finite instances of $\Theory{\CPG}$ to terms of $\DWD$.
\end{ex}

\begin{defn}\label{def:cpg_operad}
The \emph{operad of circular port graphs}, denoted $\CPG$, is the suboperad of $\DWD$ whose types are pairs of the form $\lens{X_\mathsf{port}}{X_\mathsf{port}}$ and whose terms are generated by finite instances of $\Theory{\CPG}$. 
\end{defn}

Circular port graphs differ from directed wiring diagram only in that they do not distinguish between in-ports and out-ports. This restriction is practical since circular port graphs are captured by a simpler data structure, and they formalize the composition syntax used in stencil-based numerical algorithms. We prove that every term of $\CAT{CPG}$ can be represented by a finite instance of $\Theory{\CAT{CPG}}$ in~\cite[Proposition A.3]{libkind2021operadic}.

\section{Algebras for Composing Dynamical Systems}\label{sec:composing}

Scientists often use diagrams informally to represent relationships between the components of a system. Over the last decade, applied category theorists have formalized these notions of compositional and hierarchical dynamical systems. The categorical frameworks offer a methodology for scientific modeling and their categorical structures can be implemented as modeling tools. Techniques for formalizing the composition of open dynamical systems often follow a general strategy. (1) An operad captures the syntax of interacting systems, and operadic substitution nests syntactic terms to give a more fine-grained description of the interactions. (2) An algebra over the operad assigns a concrete interpretation to the syntactic diagrams. To each type, the algebra gives a set of models of that type. To each term, the algebra gives a function that defines how to compose the chosen models.

In this section, we define operads and operad algebras for the syntax and semantics of composing open dynamical systems. The  algebras are denoted $\Dynam{\text{sem}}{\text{syn}}$ where $\text{sem} \in \{\mathsf{D}, \mathsf{C}\}$ indicates the model semantics (discrete or continuous) and $\text{syn} \in \{\to, \multimap\}$ indicates the composition syntax (directed or undirected). The algebras are all defined by lax monoidal functors, and we use the same notation for both the lax monoidal functor and the underlying algebra. For the proofs of propositions in this section, see~\cite[Appendix A]{libkind2021operadic}.

\subsection{Directed Composition} \label{sec:composition_directed}

The framework for directed composition of dynamical systems relies heavily on the generalized lens construction defined in~\cite{spivak2020GeneralizedLens}. Although this theory is robustly developed in~\cite{jazmyers2021DoubleCategories}, here we present a variant that aligns with our Julia implementation. In particular, we restrict our attention to dynamical systems defined on Euclidean spaces. 

In this section, we define algebras over $\DWD$, the operad underlying $\Lens{\Cospan(\Finset\op)}$ described in Example~\ref{ex:operad-dwd}, which factor through the algebras $\mathsf{DS}$ and $\mathsf{CS}$ defined in~\cite{spivak2016steadystates}. We define the algebras in two steps. First, let $\Euc$ be the full subcategory of the category of smooth manifolds spanned by the Euclidean spaces. A strong monoidal functor $\ev_{\Rb}: \Cospan(\Arity) \to \Euc$ is defined on objects by $P \mapsto \Rb^P$ and on morphisms by $P \xleftarrow{f} W \xrightarrow{g} Q$ maps to $g_* \circ f^*: \Rb^P \to \Rb^Q$ ~\cite[Proposition A.5]{libkind2021operadic}. By the  functoriality of the $\mathsf{Lens}$ construction, $\ev_{\Rb}$ induces a strong monoidal functor $\Lens{\ev_{\Rb}}: \Lens{\Cospan(\Arity)} \to \Lens{\Euc}$. Next, consider the lax monoidal functor $$\int_{S: \Finset} \Lens{\Euc}\left(\lens{\Rb^S}{\Rb^S}, -\right): \Lens{\Euc} \to \Set.$$ Explicitly, this functor maps an object $\lens{\Rb^{X_\tin}}{\Rb^{X_\tout}}$ to the set of pairs $\left(S \in \Finset, \lens{u}{r}: \lens{\Rb^S}{\Rb^S} \leftrightarrows \lens{\Rb^{X_\tin}}{\Rb^{X_\tout}}\right)$ and maps a morphism  $\lens{f^\#}{f}$ to the set map sending $(S, \lens{u}{r})$ to $(S,\lens{f^\#}{f} \circ \lens{u}{r})$. Finally, consider the composite 
\begin{equation}\label{eq:dwd_algebra}
    \Lens{\Cospan(\Arity)} \xrightarrow{\Lens{\ev_{\Rb}}} \Lens{\Euc} \xrightarrow{\int_{S: \Finset} \Lens{\Euc}\left(\lens{\Rb^S}{\Rb^S}, -\right)}\Set.
\end{equation}

This $\DWD$ algebra maps an object $\lens{X_\tin}{X_\tout}$  to the set of pairs $\left(S \in \Finset, \lens{u}{r}: \lens{\Rb^S}{\Rb^S} \leftrightarrows \lens{\Rb^{X_\tin}}{\Rb^{X_\tout}}\right)$. We interpret $S$ as a set of state variables and $\Rb^S$ as the state space. Depending on whether we interpret $u: \Rb^{X_\tin}  \times \Rb^S  \to \Rb^{S}$ as an indexed endomorphism of the state space or as an indexed vector field on the state space, the algebra represents either  discrete dynamical systems or  continuous dynamical systems. In other words, for an input $a \in \Rb^{X_\tin}$ and state $x \in \Rb^S$, we can either think of $u(a,x) \in \Rb^S$ as the next state or as defining the vector $\dot x = u(a,x)$. We denote the algebra defined by the composite in Equation~\ref{eq:dwd_algebra} by either $\Dynam{D}{\to}$  or $\Dynam{C}{\to}$  to highlight these distinct interpretations.\footnote{In the literature, the algebra of continuous dynamical systems explicitly represents a vector field as a section $u: \Rb^S \to T\Rb^S$ of the tangent bundle. However,  we present $\Dynam{C}{\to}$ and $\Dynam{D}{\to}$ by the same algebra because they are implemented identically in AlgebraicDynamics.} 

\subsection{Undirected Composition}

Just as there is an algebra of continuous systems over the directed syntax $\DWD$, there is an algebra over the undirected syntax $\UWD$, the operad underlying $\Cospan(\Finset)$ defined in Example~\ref{ex:operad-uwd}. We define the algebra $\Dynam{C}{\multimap}: \Cospan(\Finset) \to \Set$  which on objects takes $M$ to the set of triples $(S \in \Finset, v: \Rb^S \to \Rb^S, p: M \to S)$ and on morphisms maps the cospan $f = M \xrightarrow{q} R \xleftarrow{r} N$ to the set map $\Dynam{C}{\multimap}(f): \Dynam{C}{\multimap}(M)  \to  \Dynam{C}{\multimap}(N)$ given by $$\Dynam{C}{\multimap}(f)(S, v, p )= (S +_M R, \tilde q_* \circ v \circ \tilde q^*, \tilde p \circ r)$$ where $\tilde q$ and $\tilde p$ are defined by the pushout:
\vspace{-0.45in}

\[\begin{tikzcd}[row sep = 10, column sep = 10]
	&& {} \\
	{} & M && N \\
	S & {} & R \\
	& {S+_M R}
	\arrow["p"', from=2-2, to=3-1]
	\arrow["q", from=2-2, to=3-3]
	\arrow["r"', from=2-4, to=3-3]
	\arrow["{\tilde q}"', dashed, from=3-1, to=4-2]
	\arrow["{\tilde p}", dashed, from=3-3, to=4-2]
	\arrow["\lrcorner"{anchor=center, pos=0.125, rotate=135}, draw=none, from=4-2, to=3-2]
\end{tikzcd}\]

\noindent This composition of continuous systems is an operadic perspective of the hypergraph category $\mathsf{Dynam}$ presented in~\cite{baez2017CompositionalFramework}. Next, we define undirected composition of discrete systems. 

\begin{restatable}[]{prop}{dynamDiscreteUndirected}
There is an algebra $\Dynam{D}{\multimap}: \Cospan(\Finset) \to \Set$  which  on objects maps $M$ to the set of triples $(S \in \Finset, u: \Rb^S \to \Rb^S, p: M \to S)$ and on morphisms maps the cospan $f = M \xrightarrow{q} R \xleftarrow{r}N$ to the set map $\Dynam{D}{\multimap}(f): \Dynam{D}{\multimap}(M) \to \Dynam{D}{\multimap}(N)$ defined by $$\Dynam{D}{\multimap}(f)(S, u, p) = (S+_M R, 1_{\Rb^{S+_M R}}  + \tilde q_* \circ (u - 1_{\Rb^S}) \circ \tilde q^*  , \tilde p \circ r).$$ 
\end{restatable}

\subsection{Functorial Analysis}\label{sec:func_analysis}
Compositional modeling paves the way for compositional analysis.  Informally, an analysis of an  algebra $F: \Oa \to \Set$ is an algebra $G: \Oa \to \Set$ and a natural transformation $\blacksquare: F \Rightarrow G$ which  obscures the details of the system and highlights some feature of the behavior. Examples include identifying fixed points and orbits, solving trajectories, and computing approximations \cite{baez2017CompositionalFramework, spivak2017PixelArrays, deville2019RungeKuttaNetworks, ngotiaoco2017CompositionalityRungeKutta}. The naturality of $\blacksquare$ implies that the behavior of the total system is defined by the behaviors of its components.

For both undirected and directed dynamical systems, there exists a natural transformation which performs Euler's method. 
For a map $v: \Rb^M \times \Rb^S \to \Rb^S$ and step size $h \in \Rb_+$, define the map $\Eul_h(v): \Rb^M \times \Rb^S \to \Rb^S$ by $\Eul_h(v)(u_0, x_0) = x_0 + h v(u_0, x_0)$.

\begin{prop}[Euler's method for directed systems \cite{spivak2016steadystates}]
For $h \in \Rb_+$, there exists a natural transformation $\Eul_h^{\to}: \Dynam{C}{\to} \Rightarrow \Dynam{D}{\to}$ with components $\Eul_h^{\to}\lens{X_\tin}{X_\tout}: \Dynam{C}{\to}\lens{X_\tin}{X_\tout} \to \Dynam{D}{\to}\lens{X_\tin}{X_\tout}$  defined by $$\Eul_h^{\to}\lens{X_\tin}{X_\tout}\left(S ,  \lens{v}{r}: \lens{\Rb^S}{\Rb^S} \leftrightarrows \lens{\Rb^{X_\tin}}{\Rb^{X_\tout}}\right) = \left(S, \lens{\Eul_h(v)}{r}\right).$$
\end{prop}

\begin{restatable}[Euler's method for undirected systems]{prop}{eulerUndirected}
For $h \in \Rb_+$, there exists a natural transformation $\Eul_h^\multimap: \Dynam{C}{\multimap} \Rightarrow \Dynam{D}{\multimap}$ with components $\Eul_h^\multimap(M): \Dynam{C}{\multimap}(M) \to \Dynam{D}{\multimap}(M)$  defined by $\Eul_h^\multimap(M)(S, v, p) = (S, \Eul_h(v), p)$.
\end{restatable}

\section{Julia Implementation}\label{sec:implementation}

\begin{table}
    \centering
    \begin{tabular}{|l|l | l|}\hline 
        \textbf{Terminology} &\textbf{Mathematical abstraction} & \textbf{Julia implementation} \\ \hline
        diagram of systems & $\phi \in \Oa(s_1,\dots, s_n; t)$ & \mintinline{julia}{diagram::ACSet{TheoryO}}\\
        &$\phi_\textrm{inner} \in \Oa(r_1,\dots, r_m; s_i)$ & \mintinline{julia}{inner_diagram::ACSet{TheoryO}}\\ \hline 
        elementary models &$(m_1,\dots, m_n) \in Fs_1 \times\cdots \times Fs_n$ & \mintinline{julia}{models::Vector{T}} \\ \hline
        composition & $F(\phi)(m_1,\dots, m_n) \in Ft$ & \mintinline{julia}{oapply(diagram, models)} \\
        of models & & \\\hline
        hierarchical&$\phi \circ_i \phi_\textrm{inner}$ & \mintinline{julia}{ocompose(diagram,i,inner_diagram)} \\ 
        diagram &&\\\hline
    \end{tabular} 
    
    \caption{Comparing mathematical abstractions with their Julia implementation. Note that \mintinline{julia}{diagram}  stores the number of domain types $n$ and the types $s_1,\dots, s_n, t$. Likewise \mintinline{julia}{models}  stores the types $s_1,\dots, s_n$. The definitions of \mintinline{julia}{ocompose} and \mintinline{julia}{oapply} check that the arguments have appropriate types.}
    \label{fig:math_julia_table}
\end{table}

 In Julia, the specification of an algebra $F: \Oa \to \Set$ consists of

\begin{itemize}
    \item a schema \mintinline{julia}{TheoryO} such that finite instances of the schema represent terms of $\Oa$,
    \item a method \texttt{ocompose} implementing operadic substitution,
    \item a Julia type \texttt{T} such that values of type \texttt{T} implement algebra elements $(t \in \ob \Oa, m \in Ft)$,
    \item and a method \texttt{oapply} implementing the action of  $F$ on terms.
\end{itemize}

We highlight the correspondence between the mathematical and modeling terminologies. Let $\phi: s_1,\dots, s_n \to t$ be a term in $\Oa$. We say that $\phi$ represents a diagram or composite of subsystems. Let $m_i \in F(s_i)$ for $i = 1,\dots, n$.  The elements $(s_i, m_i)$ of $F$ are called models. We often refer to them as component models or elementary models to emphasize their role in the expression $m = F(\phi)(m_1,\dots, m_n) \in Ft$. Likewise, we say that $m$ is the composite model or the total model. The correspondences between the modeling terminology, mathematical abstractions, and Julia code constructs are listed in Table~\ref{fig:math_julia_table}.

The AlgebraicJulia ecosystem is a family of tools built on categorical techniques. Catlab.jl  provides the schemata defined in Examples~\ref{ex:theory-uwd}, ~\ref{ex:theory-dwd}, and ~\ref{ex:theory-cpg}, along with corresponding \mintinline{julia}{ocompose} methods.  AlgebraicDynamics.jl provides Julia types and \mintinline{julia}{oapply} methods that implement the models and their composition for the $\mathsf{Dynam}$ algebras of Section~\ref{sec:composing}. The following examples illustrate how applied category theory provides flexible, powerful abstractions for implementing scientific models.

\begin{ex}\label{ex:dwd-sir}

\begin{figure}[htb]
    \centering
    \begin{subfigure}[c]{0.54\textwidth}
        \includegraphics[width=\textwidth]{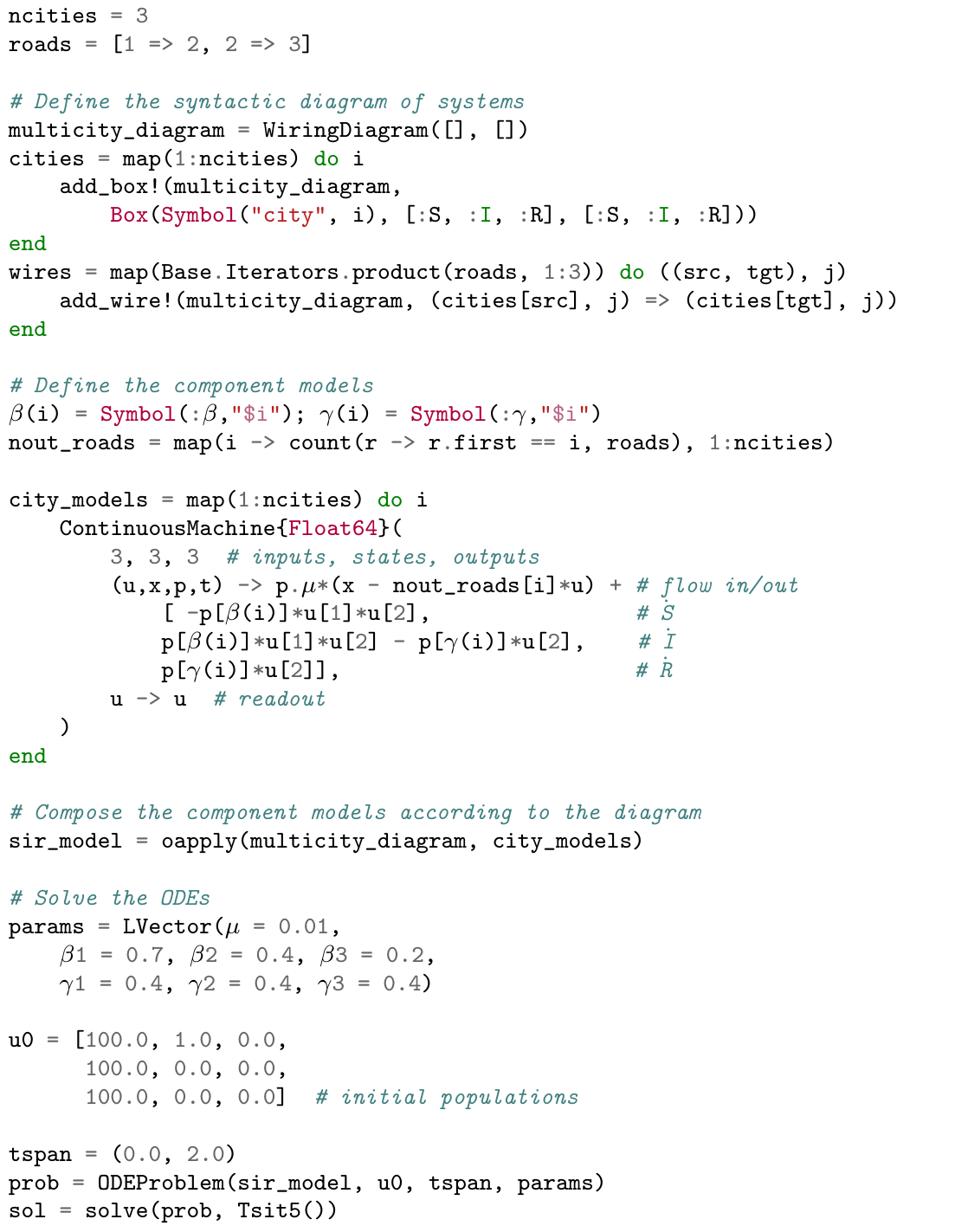}
        \subcaption[]{}
    \end{subfigure}
    \begin{subfigure}[c]{0.45\textwidth}
        \begin{subfigure}[t]{\textwidth}
            \includegraphics[width=\textwidth]{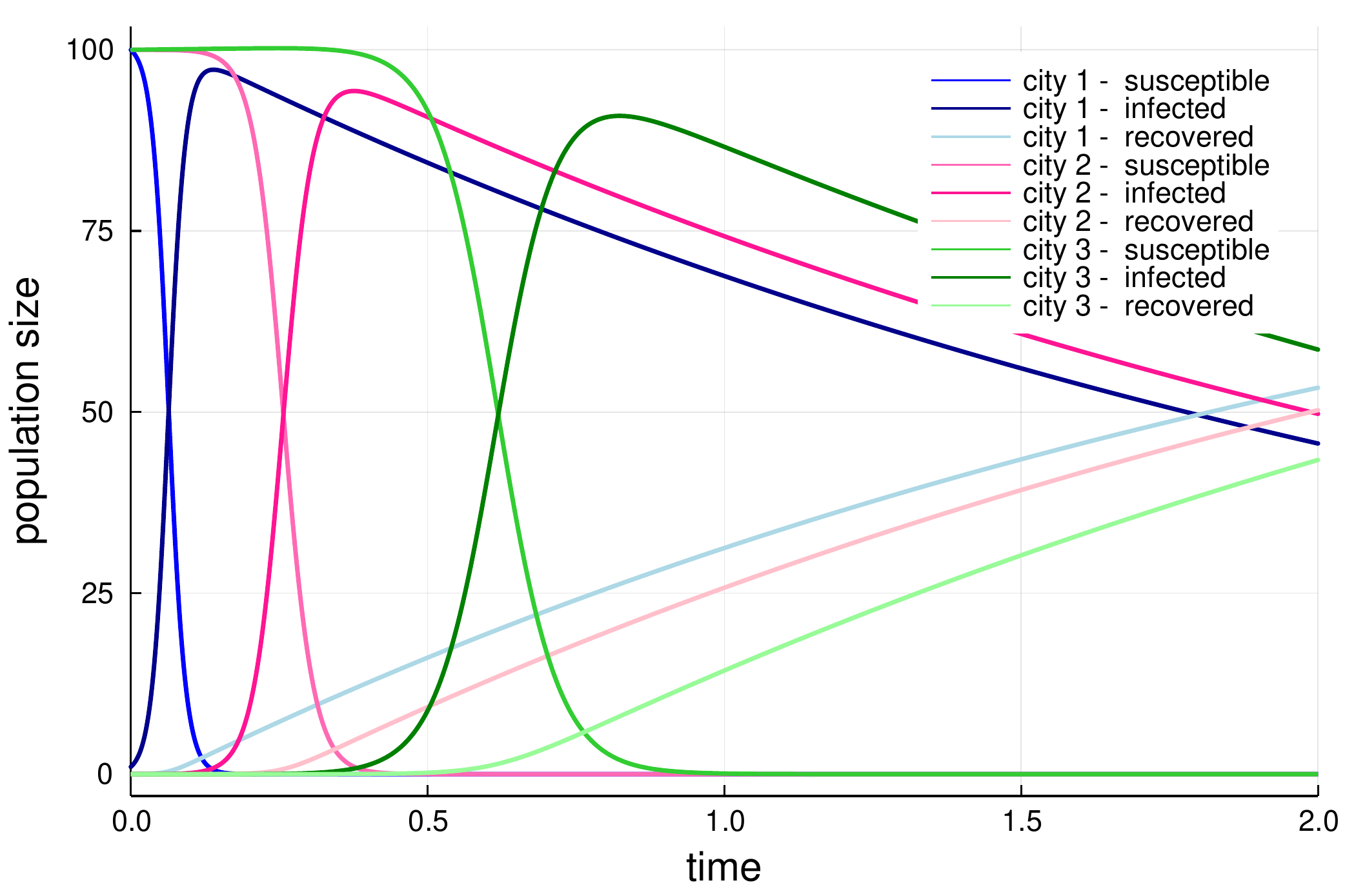}
            \subcaption[]{}
        \end{subfigure}
        
        \vspace{1cm}
        
        \begin{subfigure}[c]{\textwidth}
            \includegraphics[width=\textwidth]{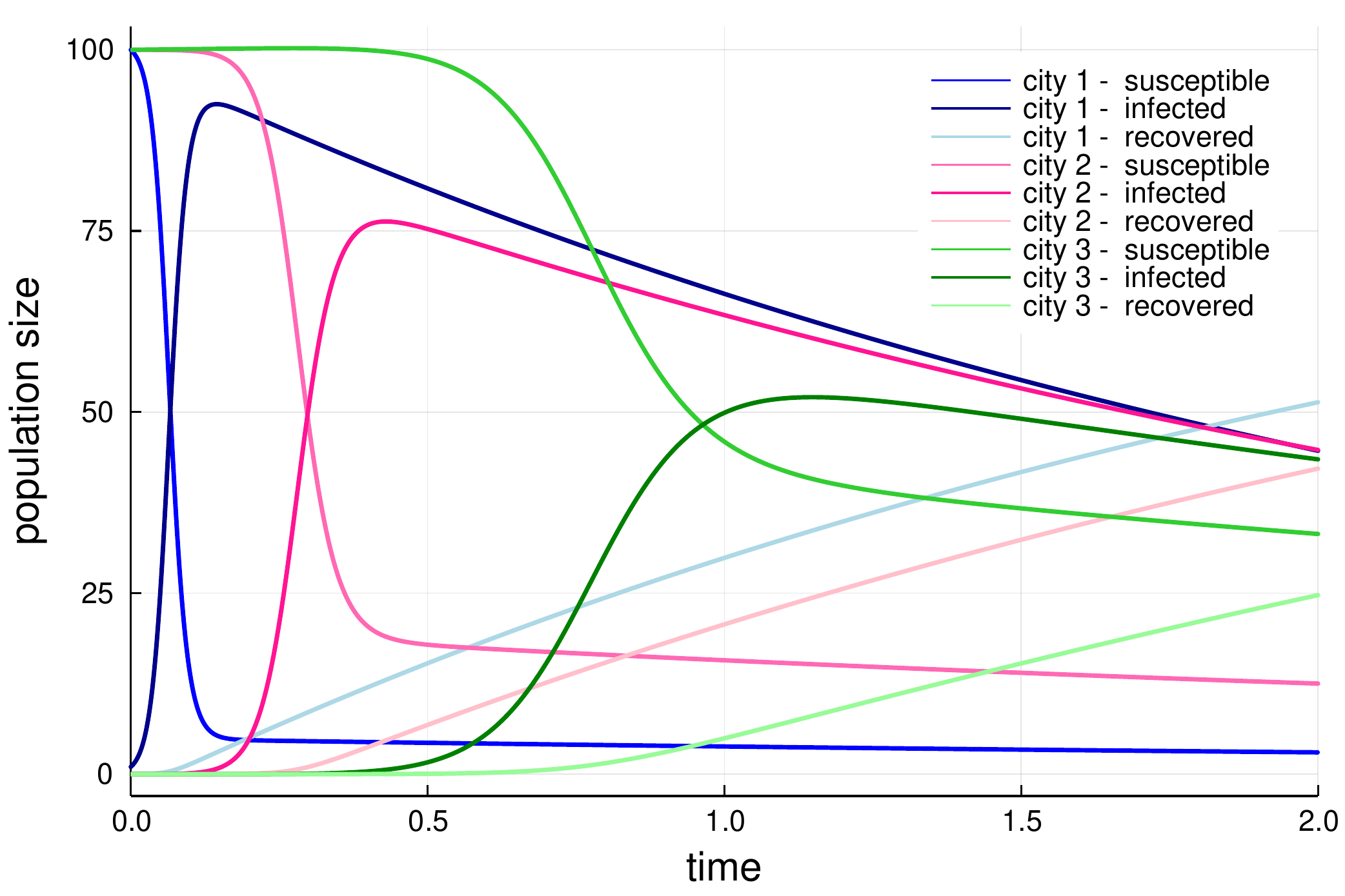}
            \subcaption[]{}
        \end{subfigure}
        
    \end{subfigure}
    
    \caption{(a) Julia code for a multi-city SIR model using AlgebraicDynamics. (b) The solution to the multi-city SIR model defined in (a). (c) The solution to a multi-city SIR-Q model produced by changing \mintinline{julia}{city_models} in the code in (a) from local SIR models to local SIR-Q models while leaving the diagram of systems intact. Comparing the plots in (b) and (c)  indicates that quarantining reduces the peak infection size in each city.}
    \label{fig:sir_code}
\end{figure}

\begin{figure}[htb]  
    \centering
    \begin{subfigure}[b]{0.15\textwidth}
        \includegraphics[scale=.6]{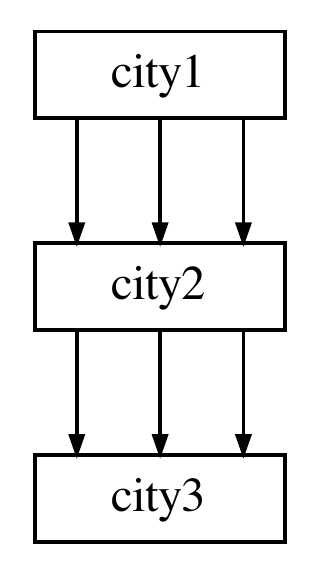}
        \subcaption[]{}
    \end{subfigure}
    \begin{subfigure}[b]{0.3\textwidth}
        \includegraphics[scale=.6]{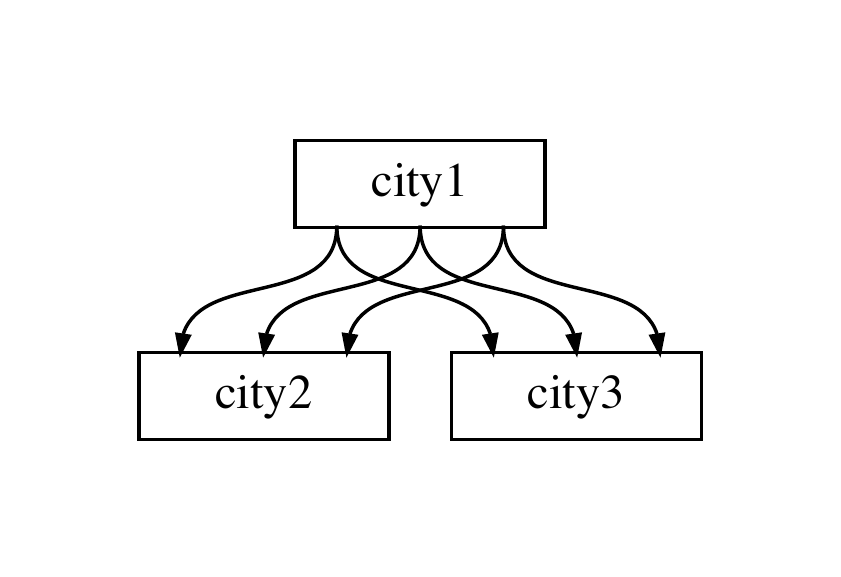}
        \subcaption[]{}
    \end{subfigure}
    \begin{subfigure}[b]{0.3\textwidth}
        \includegraphics[scale=.6]{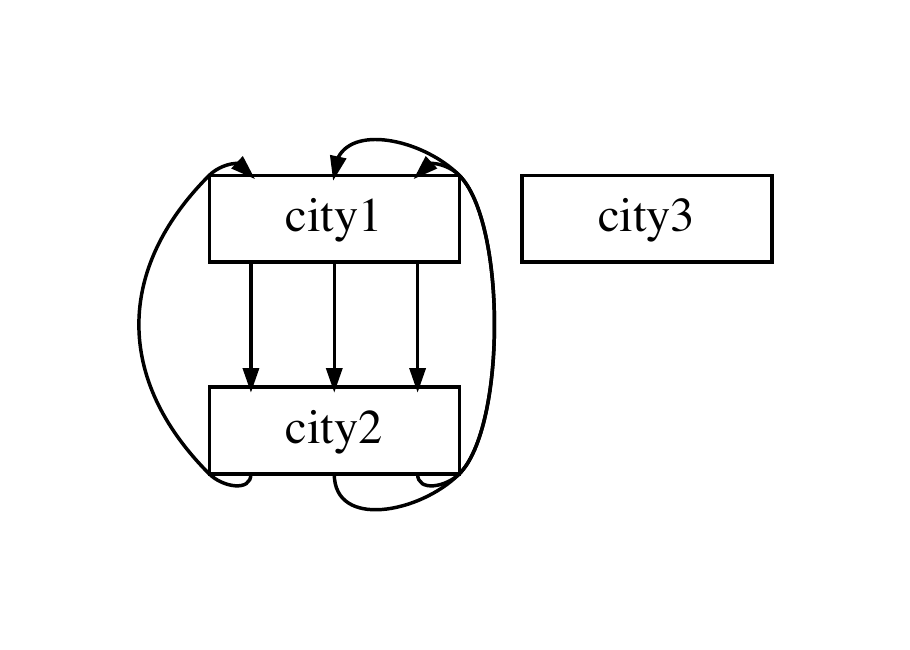}
        \subcaption[]{}
    \end{subfigure}
    \begin{subfigure}[b]{0.2\textwidth}
        \includegraphics[scale=.6]{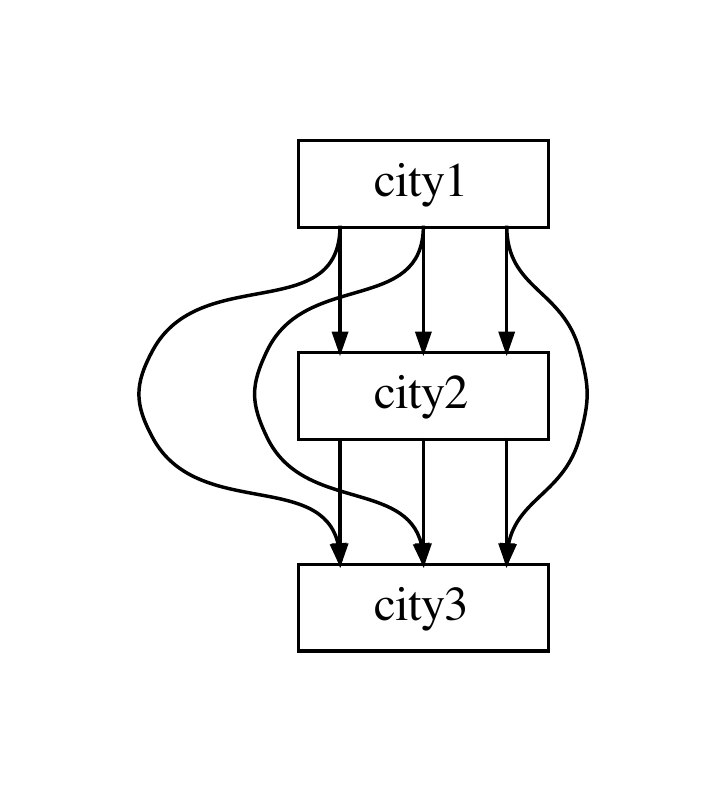}
        \subcaption[]{}
    \end{subfigure}

    \caption{The graphical display of four possible terms for composition for a multi-city SIR model with three cities. These diagrams were automatically produced by Catlab.jl. In the code in Figure~\ref{fig:sir_code} (a),  \mintinline{julia}{multicity_diagram} is defined to be the term depicted in (a).}
    \label{fig:city_diagrams}
\end{figure}

\begin{figure}
    \centering
    \begin{subfigure}[b]{0.32\textwidth}
        \includegraphics[width=\textwidth]{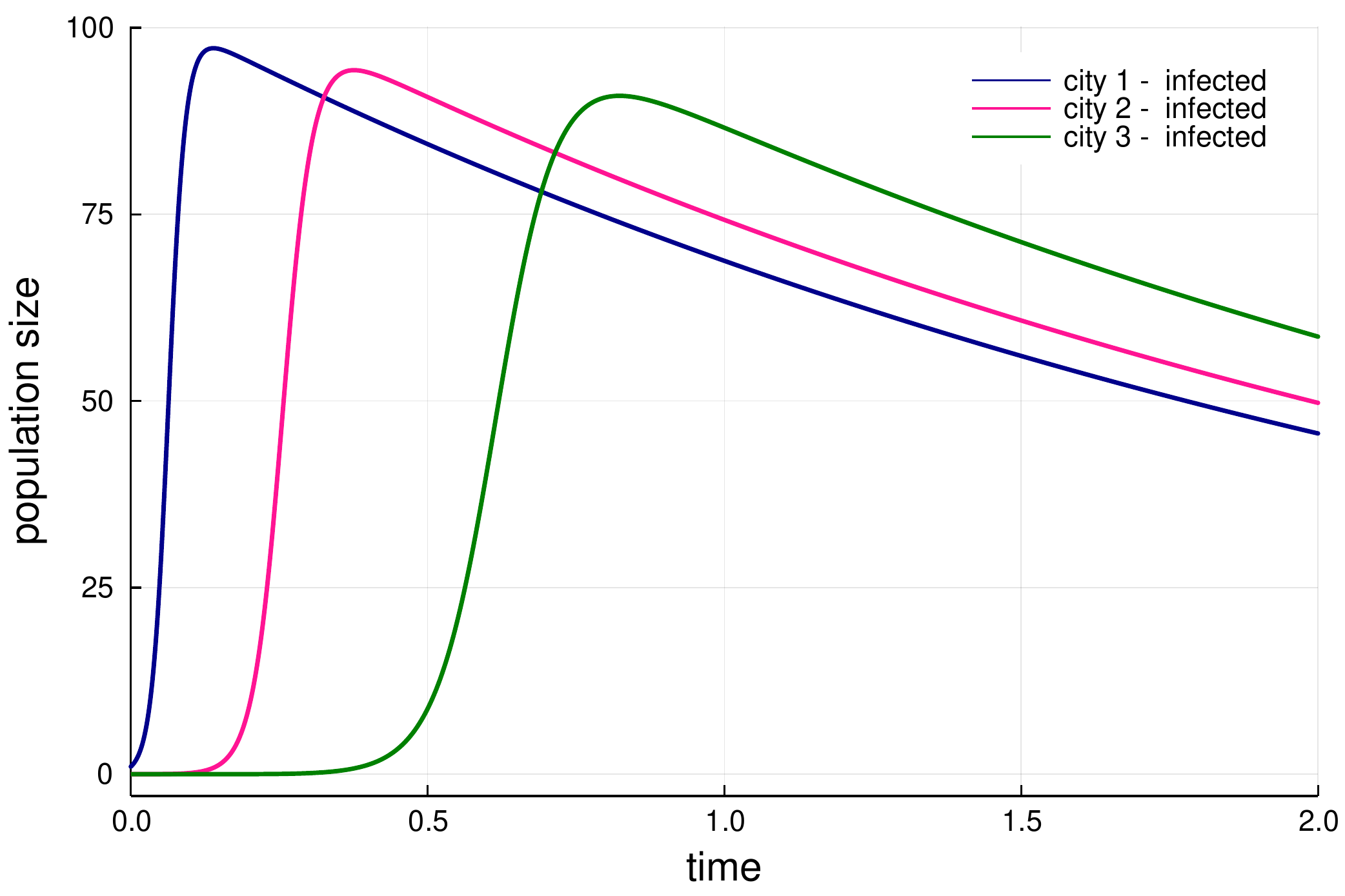}
        \subcaption[]{}
    \end{subfigure}
    \begin{subfigure}[b]{0.32\textwidth}
        \includegraphics[width=\textwidth]{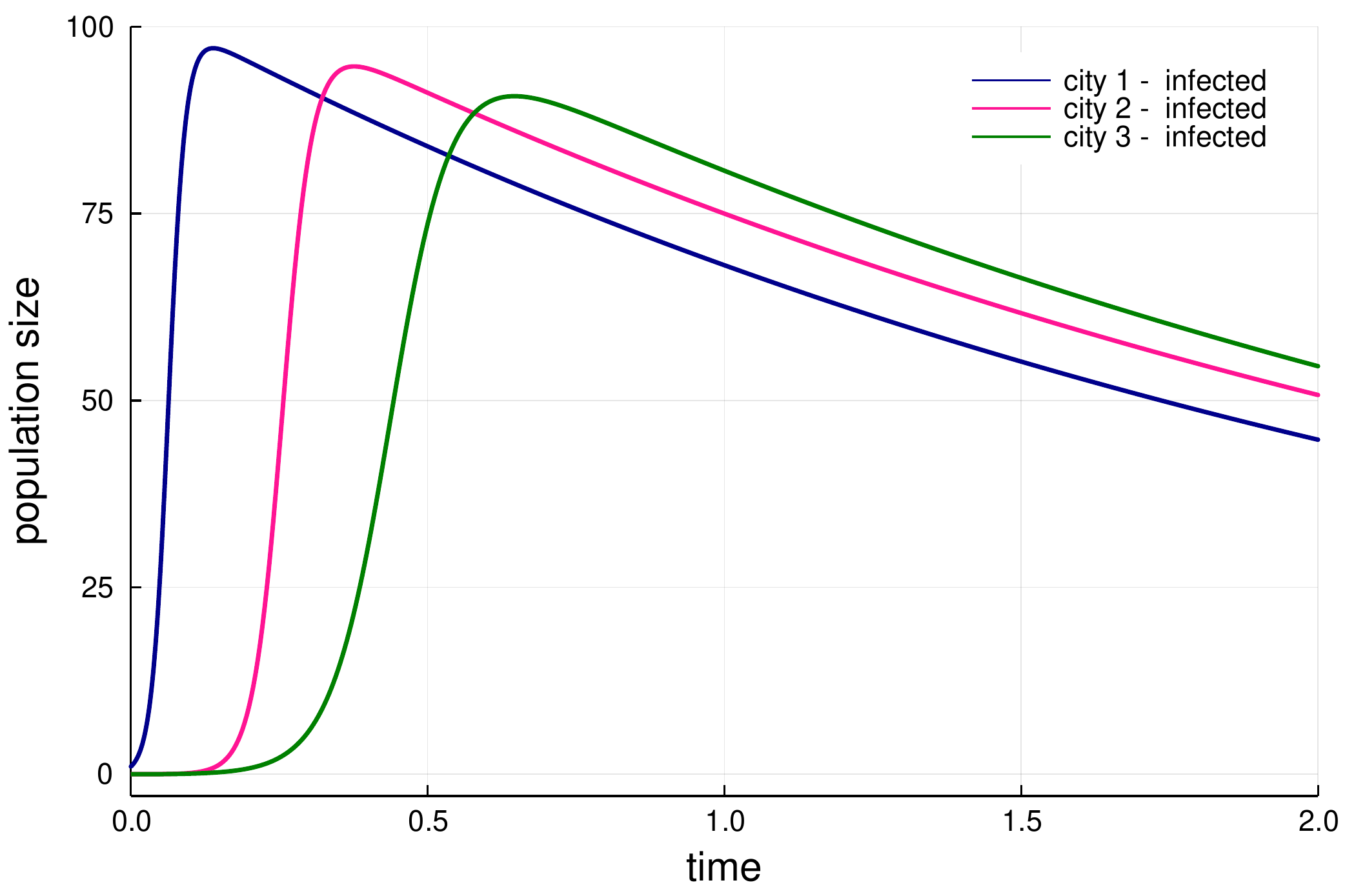}
        \subcaption[]{}
    \end{subfigure}
    \begin{subfigure}[b]{0.32\textwidth}
        \includegraphics[width=\textwidth]{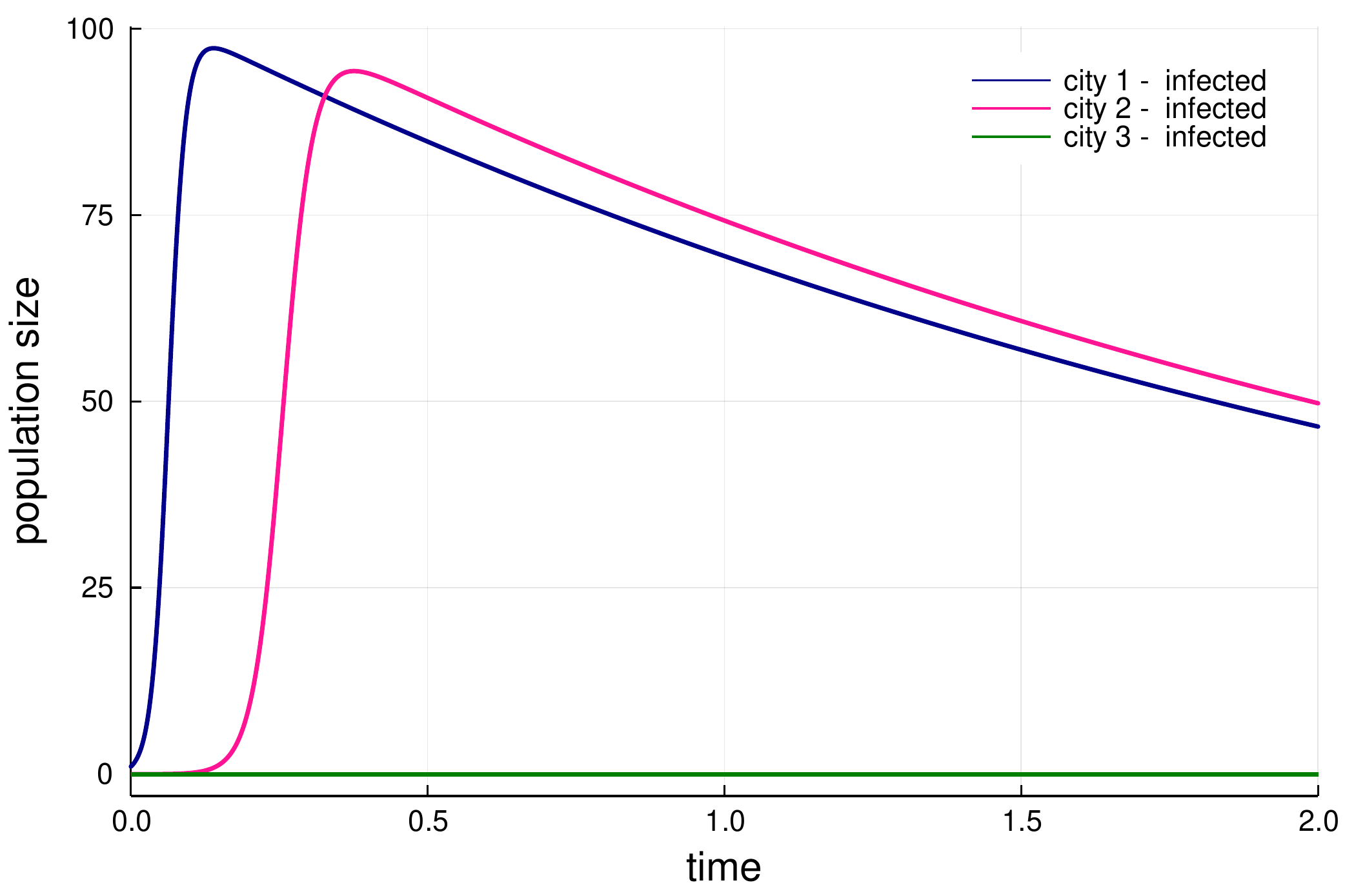}
        \subcaption[]{}
    \end{subfigure}
    \caption{Plots of each city's infected population over time for simulations of multi-city SIR models constructed using the diagrams of systems shown in Figure~\ref{fig:city_diagrams}. The plot in (a) corresponds to the multi-city SIR model defined by the diagram of systems in Figure~\ref{fig:city_diagrams}(a). Likewise for (b) and (c). These solutions reflect features of the diagrams of systems. For example, the peak in the infected population in city 3 occurs earlier in (b) than in (a) because in the diagram depicted in Figure~\ref{fig:city_diagrams}(b) city 3 is receiving infected people directly from city 1 instead of mediated by city 2.  Second, in (c) city 3 has a constant infected population of 0 because it is isolated from the epidemic occurring in cities 1 and 2. }
    \label{fig:city_solutions}
\end{figure}

The SIR model is a classic model of the spread of an infectious disease \cite{kermack1927contributionmathematical}. A single-city SIR model has a susceptible population $S$, an infected population $I$, and a recovered population $R$, which evolve according to the continuous dynamics 
\begin{equation}\label{eq:sir}
    \dot S = -\beta SI, \qquad \dot I = \beta SI - \gamma I, \qquad \dot R = \gamma I.
\end{equation} Such a model assumes a mixing of the population which fails to account for geographic or social distinctions between sub-populations. To model such distinctions, we  compose multiple single-city SIR models using the operad algebra $\Dynam{C}{\to}: \DWD \to \Set$ defined in Section~\ref{sec:composition_directed}. 

Figure~\ref{fig:sir_code}(a) shows the code for a multi-city SIR model using AlgebraicDynamics. The directed wiring diagram, \mintinline{julia}{multicity_diagram}, is an instance of $\Theory{\CAT{DWD}}$ and so corresponds to a term $\phi$ in $\DWD$. Four different wiring diagrams for composing three single-city SIR models are shown in Figure~\ref{fig:city_diagrams}, of which Figure~\ref{fig:city_diagrams}(a) is used in the code in Figure~\ref{fig:sir_code}(a). The Julia type \mintinline{julia}{ContinuousMachine} implements elements of $\Dynam{C}{\to}$, and the array \mintinline{julia}{city_models} defines for each city a model with three states (corresponding to the local $S$, $I$, and $R$ populations), the vector field from Equation~\ref{eq:sir} modulated by an inflow and outflow, and an identity readout. The call of the function \mintinline{julia}{oapply} applies the set map $\Dynam{C}{\to}(\phi)$ to the single-city elements defined by \mintinline{julia}{city_models} and returns a multi-city model \mintinline{julia}{sir_model}. Figure~\ref{fig:sir_code}(b) shows a solution to the multi-city model. 

Several aspects of the algebraic formalism translate to practical software features. First, the clear delineation between syntax and semantics enables the user to modify them independently. For example, the choice of the composition term \mintinline{julia}{multicity_diagram} is independent of the choice of the single-city models \mintinline{julia}{city_models}. Figure~\ref{fig:sir_code}(c) gives an example of leaving the syntax intact while modifying the semantics from local SIR models to local SIR-Q models which represent quarantining populations in each city. Conversely, Figure~\ref{fig:city_solutions} gives three examples of modifying the syntactic diagram of cities while leaving the local SIR model intact. Second, adding more components, such as additional cities, to an existing model is a straightforward application of compositionality. Without the algebraic abstraction, this procedure would involve many coordinated changes to the code, an error-prone and time-consuming process. Finally, the syntactic diagram gives information about the behavior of the model. For example, the composition term in Figure~\ref{fig:city_diagrams}(c) implies that the behavior of city 3 is independent of cities 1 and 2 regardless of the choice of component models. Therefore, any analysis of these subsystems can be done in parallel. 
\end{ex}

\newpage
\begin{ex}\label{ex:uwd-eco}
A model of an ecosystem is composed of many primitive growth, decline, and predation models where a single species may be involved in multiple subsystems. The compositional approach allows us to divide the problem into interacting subsystems and conquer them independently. At the highest level there are two components to the total ecosystem: a land system and a river system. These subsystems compose by identifying species which appear in both. In this case, the hawk populations belonging to the land and river systems are identified in the total system. Figure~\ref{fig:eco_diagrams}(a) depicts the corresponding term in the operad $\UWD$.

\begin{figure}[bht]
\footnotesize\centering
\begin{minipage}[b]{0.5\textwidth}
\centering
\includegraphics[width=\textwidth]{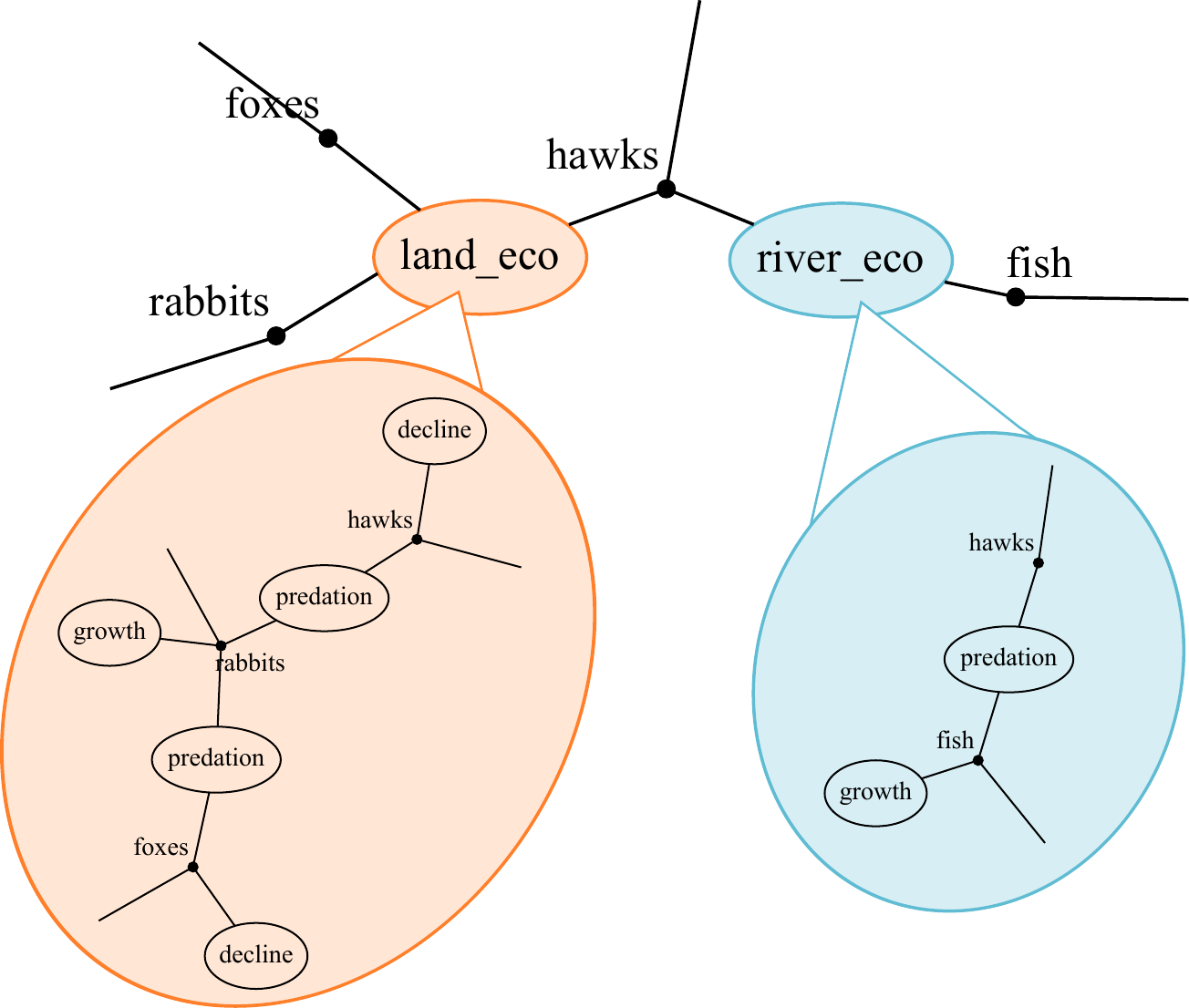}
(a)
\end{minipage}\hfill
\begin{minipage}[b]{0.35\textwidth}
\centering
\includegraphics[width=\textwidth]{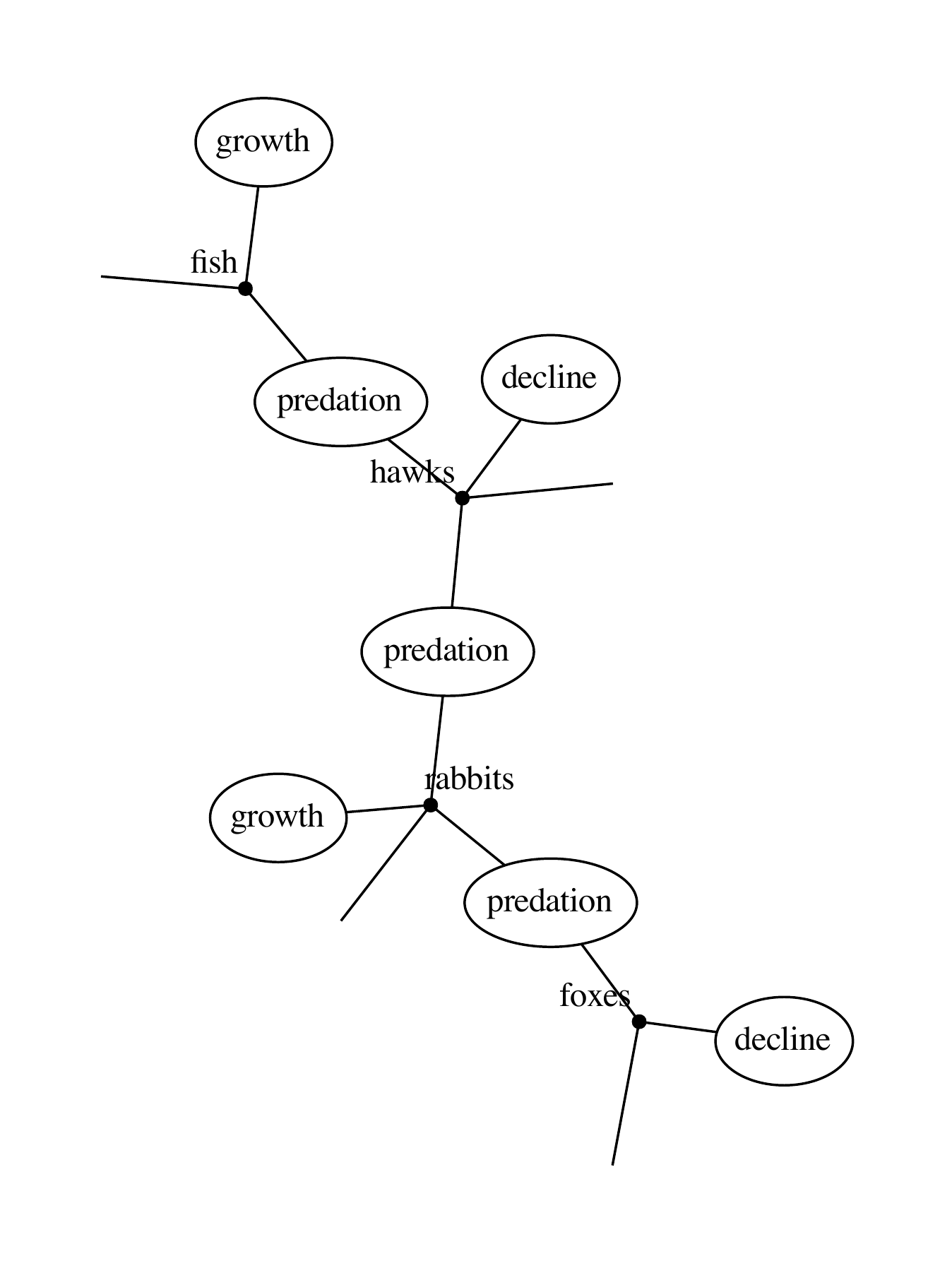}
(b)
\end{minipage}

    \caption{(a) The graphical depiction of three terms of $\UWD$ and their hierarchical relationships. Top: \mintinline{julia}{total_diagram}. Bottom left: \mintinline{julia}{land_diagram}. Bottom right: \mintinline{julia}{river_diagram}. The graphics (in black) are produced by the Catlab.Graphics module. The coloring highlights the nested structure. (b) The result of substituting the terms for the land and river systems into the term for the total ecosystem is \mintinline{julia}{ocompose(total_diagram, [land_diagram, river_diagram])}.
    }
    \label{fig:eco_diagrams}
\end{figure}

\begin{figure}[htb]
    \centering
    \includegraphics[width=0.6\textwidth]{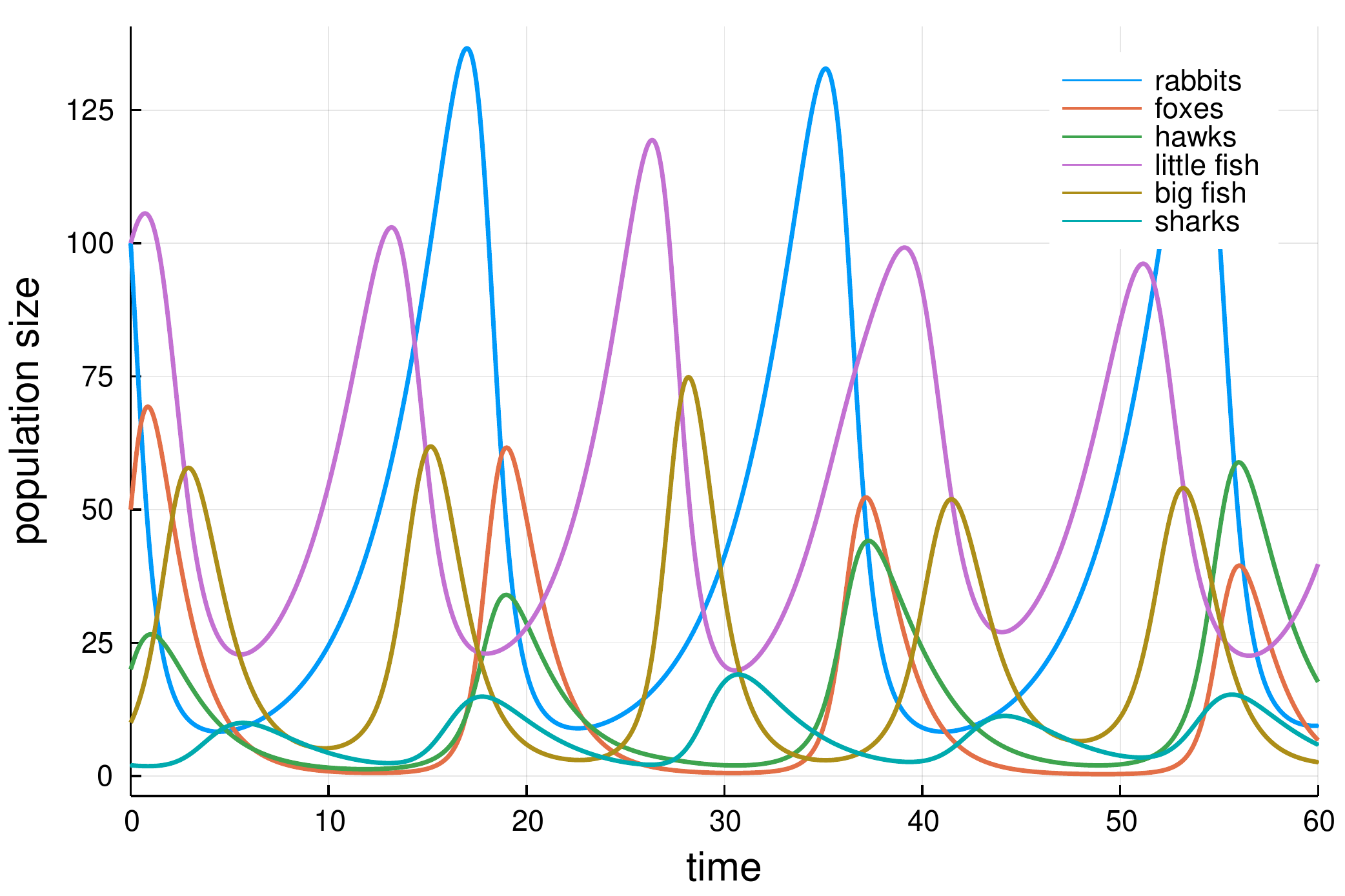}
    \caption{The solution to the complete ecosystem model defined hierarchically in Example~\ref{ex:uwd-eco}.}
    \label{fig:eco_sol}
\end{figure}

This approach to modeling is modular and hierarchical. It is modular because the models for the land and river ecosystems are defined independently. It is hierarchical because the models for the land and river ecosystems can themselves be composed of still more fine-grained interactions. Figure~\ref{fig:eco_diagrams}(a) gives terms for the fine-grained interactions within the land and river subsystems respectively and shows how they nest into the high-level diagram. Furthermore, the syntax for the land subsystem can be defined as a pushout in $[\Theory{\UWD}, \Set]$ highlighting the advantages of using $\Ca$-Sets to define the composition syntax. Given this hierarchical structure, there are two different but equivalent strategies to construct the complete ecosystem model. The model elements for the growth, decline, and predation interaction types can be applied either  at the level of sub-ecosystems:
\begin{minted}[fontsize = \footnotesize]{julia}
    land_sys = oapply(land_diagram, land_models)
    river_sys = oapply(river_diagram, river_models)
    total_sys = oapply(total_diagram, [land_sys, river_sys])
\end{minted}

\noindent or  at the level of the total ecosystem:\footnote{The term of $\UWD$ corresponding to $\mintinline{julia}{eco_diagram}$ is shown in Figure~\ref{fig:eco_diagrams}(d).}

\begin{minted}[fontsize = \footnotesize]{julia}
    eco_diagram = ocompose(total_diagram, [land_diagram, river_diagram])
    total_sys = oapply(eco_diagram, vcat(land_models, river_models))
\end{minted}

Functoriality of $\Dynam{C}{\multimap}$ implies that these strategies produce models with identical denotational semantics, a solution for which is shown in Figure~\ref{fig:eco_sol}. Since syntactic terms can be built and interpreted hierarchically with equivalent results, hierarchical modeling is flexible, scalable, and parallelizable. For example, if we discover a fourth species involved in the land system, then we can adjust its internal syntax and semantics independently of the river system.

\end{ex}

\begin{figure}[bht]
\footnotesize
    \begin{minipage}{0.5\linewidth}
    \centering
    \includegraphics[width=\textwidth]{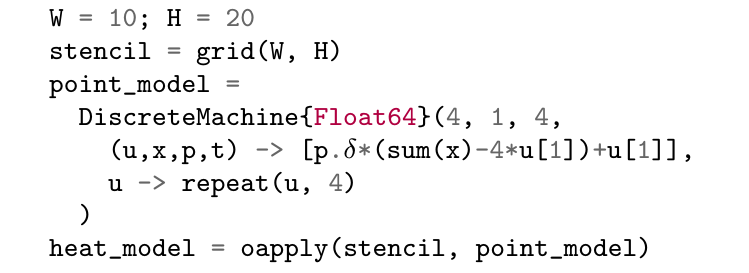}
    (a)
    \end{minipage}\begin{minipage}{0.5\linewidth}
    \centering
        (b) \includegraphics[width = 0.75\textwidth, valign=c]{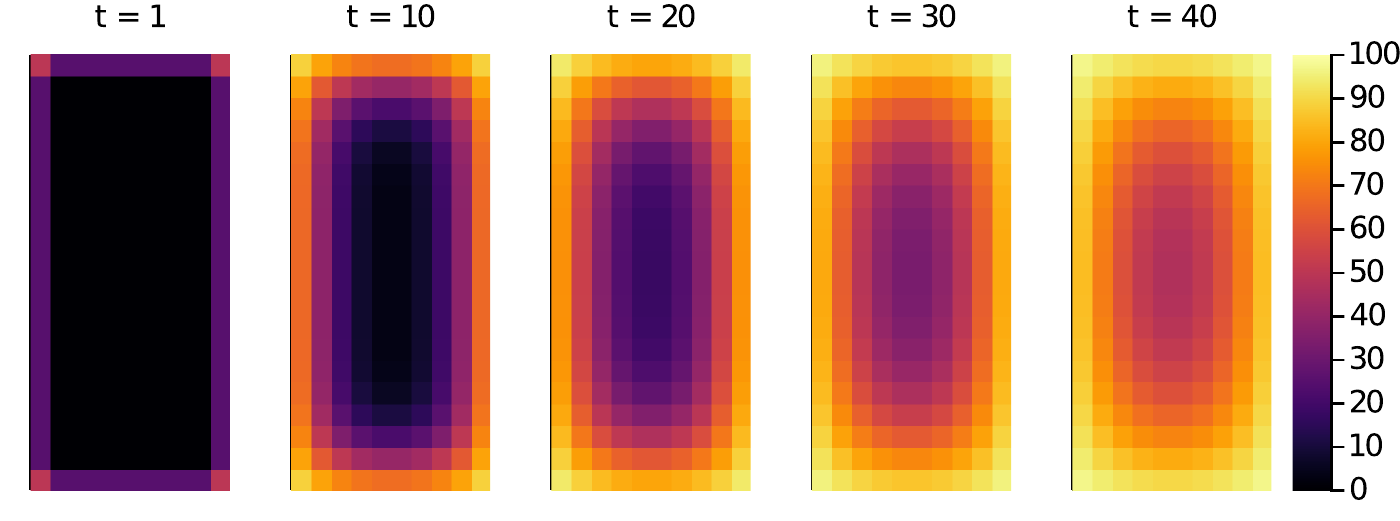}\\
        (c) \includegraphics[width = 0.75\textwidth, valign=c]{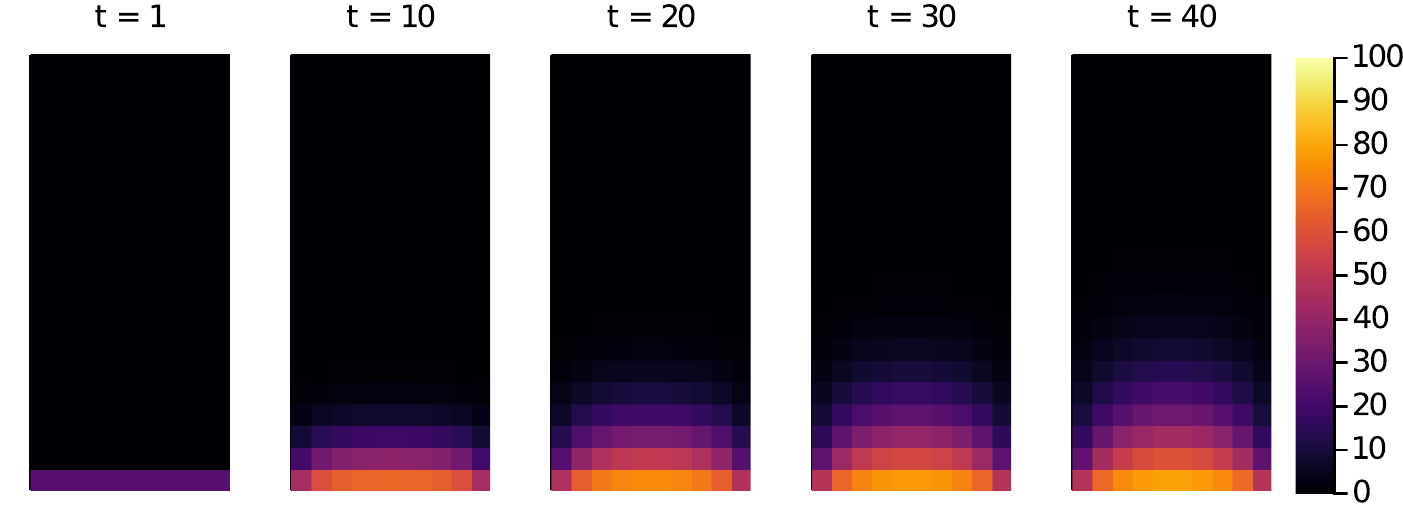}
    \end{minipage}
    \caption{(a) Julia code for a discrete approximation of the heat equation. (b) Simulation of the model over time where the boundary conditions specify that heat enters the system from all sides. (c) Simulation of the model over time where the boundary conditions specify that heat enters the system from the bottom. }
    \label{fig:heat_eq}
\end{figure}

\begin{ex}
In Definition~\ref{def:cpg_operad}, we introduced the operad of circular port graphs $\CPG$ and showed that it is a suboperad of $\DWD$. The compositions

\[\begin{tikzcd}
	\CPG & \DWD & &\Oa(\Set)
	\arrow[hook, from=1-1, to=1-2]
	\arrow["\Dynam{C}{\to}", shift left=1, from=1-2, to=1-4]
	\arrow["\Dynam{D}{\to}"', shift right=1, from=1-2, to=1-4]
\end{tikzcd}\]define algebras of continuous and discrete systems over circular port graphs. In AlgebraicDynamics, the \mintinline{julia}{oapply} method for circular port graphs is implemented independently of the \mintinline{julia}{oapply} method for directed wiring diagrams to improve performance.

The operads and operad algebras for circular port graphs   formalize standard numerical analysis techniques, namely the method of stencils and the finite difference method. Figure~\ref{fig:heat_eq} gives a solution to the 2D heat equation using the finite difference method. The method \mintinline{julia}{grid(W::Int, H::Int)}  returns a circular port graph with $W \times H$ nodes arranged in a grid. Each node has four ports and is connected to its four adjacent neighbors via symmetric wires. The unattached ports of boundary nodes are exposed by open ports. This syntactic term recovers the 5-point stencil.
\end{ex}

\begin{ex}

So far we have presented three examples of \textit{constructing} models compositionally. The next step is to \textit{analyze} models compositionally, using the mathematical abstraction of natural transformations defined in Section~\ref{sec:func_analysis}.  Let $\blacksquare: F \Rightarrow G$ be a natural transformation between operad algebras. If the Julia type \mintinline{julia}{S} implements elements of $F$ and the Julia type \mintinline{julia}{T} implements elements of $G$, then $\blacksquare$ is implemented by a Julia method \mintinline{julia}{blackbox(::S)::T}. AlgebraicDynamics.jl implements Euler's method for both undirected and directed systems.

\end{ex}

\section{Conclusion}

Applied category theory offers rigorous denotational semantics for scientific modeling. The development of Catlab and the AlgebraicJulia ecosystem provides the computational framework for implementing these semantics and making them accessible and practical for scientists. The implementation in AlgebraicDynamics of operad algebras for composing open dynamical systems  is a first example of this work and enables modelers to leverage the abstractions of applied category theory to solve real-world problems. In this paper, we showcased  many advantages that these abstractions offer to modelers, including hierarchical modeling, hierarchical model interpretation, and independence of model syntax and semantics. These advantages reduce development time for model construction, exploration, and refinement and facilitate parallel development across distinct expert domains.

In future work, we will lay a foundation for using higher category theory to study relationships between dynamical models. For example, we can use double categories of structured cospans with resource-sharing semantics to study how localized changes to model structure affect system behavior. Furthermore, while a theme of this work is ``implementations informed by abstractions," its converse ``abstractions informed by implementations" is a rich source of mathematical ideas.  For instance, we saw that Euler's method is a functorial process, yet many standard numerical methods are not functorial. We aim to develop abstractions characterizing this lossiness and use them to prove accuracy results. As another direction, the algebras presented in Section~\ref{sec:composing} formalize the concept of real-valued data flowing along wires. However, the Julia implementation suggests that these algebras generalize to types with Frobenius structure. Finally, implementing the terms of operads as $\Ca$-sets is a useful device for treating syntactic terms as data structures. We conjecture that this is an artifact of a more general theory of operads defined by $\Ca$-sets.

\bibliographystyle{eptcs}
\bibliography{ref}

\end{document}